\documentclass[proc]{edpsmath}
\usepackage{tikz}
\usetikzlibrary{calc}
\usepackage{dsfont} % Indicator function

\begin{document}

\newcommand{\E}{\mathbb{E}} % Expectation
\newcommand{\Var}{\mathbb{V} {\rm ar}} % Variance
\newcommand{\VarR}{\mathbb{V} {\rm arR}} % Relative variance

\newcommand{\eps}{\varepsilon}
\newcommand{\dps}{\displaystyle}
\newcommand{\D}{{\cal{D}}}

\newtheorem{hyp}[thrm]{\textsc{Assumption}}

%%-----------------------------
%%      the top matter
%%-----------------------------
\title{A parameter identification problem in stochastic homogenization}
\author{F. Legoll}
\address{Laboratoire Navier, \'Ecole Nationale des Ponts et Chauss\'ees, Universit\'e Paris-Est, 6 et 8 avenue Blaise Pascal, \\
77455 Marne-La-Vall\'ee Cedex 2, France;
\email{legoll@lami.enpc.fr}}
\secondaddress{INRIA Rocquencourt, MATHERIALS research-team, Domaine de Voluceau, B.P. 105,
78153 Le Chesnay Cedex, France}
\author{W. Minvielle}
\address{CERMICS, \'Ecole Nationale des Ponts et
Chauss\'ees, Universit\'e Paris-Est, 6 et 8 avenue Blaise Pascal, \\
77455 Marne-La-Vall\'ee Cedex 2, France;
\email{william.minvielle@cermics.enpc.fr}}
\sameaddress{, 2}
\author{A. Obliger}
\address{Sorbonne Universit\'es, UPMC Univ. Paris 06, UMR 8234 PHENIX, 75005 Paris, France;
\email{amael.obliger@upmc.fr}}
\secondaddress{CNRS, UMR 8234 PHENIX, 75005 Paris, France and ANDRA, Parc de la Croix-Blanche, 1-7, rue Jean-Monnet, 92298 Ch\^atenay-Malabry, France}
\author{M. Simon}
\address{UMPA, UMR-CNRS 5669, ENS Lyon, 46 all\'ee d'Italie, 69007 Lyon, France;
\email{marielle.simon@ens-lyon.fr}}
%
%\dedicated{\it Dedicated to Maurice Dupont} %if necessary
%
\begin{abstract}
In porous media physics, calibrating model parameters through experiments is a challenge. This process is plagued with errors that come from modelling, measurement and computation of the macroscopic observables through random homogenization -- the forward problem -- as well as errors coming from the parameters fitting procedure -- the inverse problem. In this work, we address these issues by considering a least-square formulation to identify parameters of the microscopic model on the basis on macroscopic observables. In particular, we discuss the selection of the macroscopic observables which we need to know in order to uniquely determine these parameters. To gain a better intuition and explore the problem without a too high computational load, we mostly focus on the one-dimensional case. We show that the Newton algorithm can be efficiently used to robustly determine optimal parameters, even if some small statistical noise is present in the system. 
\end{abstract}
\begin{resume} 
En physique des milieux poreux, calibrer certains param\`etres d'un mod\`ele microscopique sur la base d'exp\'eriences donnant acc\`es \`a des grandeurs macroscopiques est un enjeu majeur. Cette d\'emarche est entach\'ee d'erreurs de mod\`ele, de mesure et de calculs dans la proc\'edure d'homog\'en\'eisation: le probl\`eme direct est biais\'e. La r\'esolution du probl\`eme inverse, lorsqu'il s'agit d'estimer les param\`etres \`a partir des observations, engendre aussi des erreurs. Nous consid\'erons ici une formulation ``moindres carr\'es'' du probl\`eme, cherchant \`a minimiser l'erreur entre les quantit\'es macroscopiques observ\'ees et celles calcul\'ees via l'homog\'en\'eisation al\'eatoire. Nous discutons en particulier de la nature des informations macroscopiques n\'ecessaires pour d\'eterminer de mani\`ere univoque les param\`etres de la densit\'e de probabilit\'e des propri\'et\'es microscopiques. Afin d'explorer plus facilement cette question, nous nous int\'eressons ici essentiellement au cas unidimensionel. Nous montrons que le probl\`eme peut \^etre r\'esolu de mani\`ere efficace par l'algorithme de Newton, m\^eme en pr\'esence d'un petit bruit statistique. 
\end{resume}
\maketitle
%%-----------------------------
%%      your text
%%-----------------------------

\section{Introduction}

Modelling porous media is a challenge, in particular because geometry of such materials can be extremely complex. Rock samples are often described as a pile of layers 
of solid phase which do not permit flows, creating voids in-between layers that are connected by channels, the size and shape of which is difficult to describe (and to observe experimentally, although, in rare cases, imaging methods such as micro-tomography can be used). Besides these issues related to the description of the geometry of the media, another difficulty is to properly model the physical phenomena occuring in the flow. 
To circumvent these difficulties, a possible approach consists in completely forgetting the exact geometry of the system except for a few parameters (e.g. the size of the channels), and consider that the channels form a simple network, often taken to be $\xZ^d$. This results in the so-called pore-network models (PNM), initially introduced by Fatt in the 1950s~\cite{F56} and which have been widely used since then. The void space of a rock (its porosity) is described by a pore network connected by channels. In this framework, the geometry of pores and channels is idealized. Some microscopic properties are assigned to network elements (e.g. the conductance of the channels) and rules are defined to compute the upscaled (homogenized) properties on the basis of this microscopic description. In turn, these upscaled properties can be compared to the available experimental data. The aim is to construct a microscopic network with the same effective properties as those of a real representative sample of rock. 

\medskip

In this work, we follow this approach, and assume that pores are located at the vertices of a simple lattice. We adopt a stochastic model, and assume that, at the microscopic scale, physical properties are described by some random field. In particular, we focus on monophasic transport phenomena in porous media, where the sample of rock is mainly characterized by its permeability. These phenomena are described by the Darcy's law, where the local flux of water is assumed to be proportional to the local pressure gradient, and the microscopic properties of interest are the conductances of the channels. In the pore network model, conductances are solely assigned to channels, and it is assumed that pores do not contribute to the flow. Following Darcy's equation, the microscopic pressure field is computed in the network by ensuring mass conservation at each pore. The equation to solve is therefore a {\em discrete} linear elliptic equation in divergence form, with random coefficients (see~\eqref{eq:conservation}--\eqref{eq:recast} below).

The conductances of the channels (i.e. their microscopic permeabilities) depend on their size. Therefore the construction of the network starts by randomly attributing a size 
to each channel. In practice, this channel size distribution can be inferred from experiments such as mercury porosimetry: we denote it by ${\cal L}_{\rm exp}$. Several issues of different nature arise in this procedure. As a consequence, it turns out that the effective properties (e.g. macroscopic permeability) that are computed for a pore network with channel sizes distributed according to ${\cal L}_{\rm exp}$ are different from the experimental effective properties. The extraction procedure, which provides a channel size distribution, is thus somewhat slightly inconsistent. The main goal of this work consists in improving that distribution, when starting from the experimental initial guess, in order to eventually achieve a better agreement between measured and computed effective properties.

\medskip

From a more mathematical standpoint, the question can be phrased in the following terms. Consider a second-order divergence-form operator whose coefficients are random. If the distribution of the coefficients is stationary and ergodic, then (under some additional technical assumptions) this random operator can be replaced, over large scales, by an effective operator with constant homogenized coefficients. Random homogenization theory actually provides formulas to compute the homogenized quantities. We have thus at our disposal a procedure to compute macroscopic quantities if we know the microscopic quantities, and to solve the so-called {\em forward problem}. However, in practice, given a heterogeneous materials, it is a difficult question to decide on the law of the microscopic physical properties. On the other hand, macroscopic quantities are more easily accessible. It is thus of interest to consider the {\em inverse problem}, and try to extract some information on the properties of the materials at the microscopic scale on the basis of macroscopic quantities. 

In the same spirit, if one makes assumptions on the microscopic law, then macroscopic quantities can be computed, and for instance compared to experimental values. In view of the possible discrepancy between the two, one could question or revisit the assumptions made at the microscopic scale. 

Of course, homogenization is an averaging process, which filters out many features of the microscopic coefficients. There is thus no hope to recover a full information about the microstructure (in our case, the {\em probability distribution} of the conductances) from the only knowledge of macroscopic quantities. We adopt here a more restricted objective. We will assume a functional form for the distribution of the microscopic conductances (namely, a Weibull distribution). Our aim is to recover the {\em parameters} (denoted here $\theta$) of that microscopic law of the basis of macroscopic quantities.

\medskip

We point out that our approach is not specific to Weibull laws, and that it could be used for other distribution laws with parameters $\theta$. What we need is that the random field $A(x,\omega)$ used at the microscopic scale can be written as
$$
A(x,\omega) = {\cal F} \Big( u(x,\omega),\theta \Big)
$$
where $u(x,\omega)$ is a field of random variables that are uniformly distributed and ${\cal F}$ smoothly depends on the parameters $\theta$ (see~\eqref{eq:cov} in our particular case). Computing the derivatives of the microscopic random field $A(x,\omega)$ (and next of the macroscopic, homogenized quantities) with respect to $\theta$ is then easy. Our motivation for choosing Weibull laws comes from physical reasons: based on experimental results, it appears to be a reasonable choice.

Likewise, our approach is not specific to {\em discrete} elliptic equations. It could be also applied for problems modelled with {\em continuous} elliptic partial differential equations (PDEs) with random, highly oscillatory coefficients. Here, we consider discrete equations because the pore network model, which is naturally written in terms of discrete equations, is commonly used for such materials. 

\medskip

The question of recovering the unknown parameters $\theta$ of the microscopic distribution from homogenized (and more generally macroscopic) quantities belongs to the wide family of \emph{inverse problems}. In this work, a major point of interest is the selection of the macroscopic quantities which we need to know in order to uniquely determine the parameters $\theta$. This point is discussed in Section~\ref{sec:discut}.

The article is organised as follows. In Section~\ref{sec:definition}, we recall some elements of stochastic homogenization and describe the physical problem that motivates this work (including the choice of Weibull laws). We conclude that section with results specific to the one-dimensional case. In particular, random variables distributed according to a Weibull law are not isolated from 0 or $+\infty$, and thus the microstructure does not satisfy the classical assumption of ellipticity, namely~\eqref{eq:ellipticity} below. We show in Section~\ref{sec:oneD} that, in the one-dimensional case, homogenization still holds under a weaker assumption, that in turn is satisfied by Weibull random variables. 

Next, in Section~\ref{sec:inverse}, we introduce our parameter fitting problem, formulated as a least-square optimization. We first consider the general (multi-dimensional) case before turning to the one-dimensional case. In that latter case, we discuss the macroscopic quantities that are needed to uniquely determine the parameters $\theta$. More precisely, Weibull laws have two parameters, and the knowledge of a single homogenized quantity (namely the macroscopic permeability) is, as expected, insufficient to determine the two unknown parameters. We show there (in the one-dimensional case) that, if we additionally specify the relative variance of the effective macroscopic permeability, then we are in position to uniquely determine the two parameters of the microscopic Weibull law. 

Section~\ref{sec:numerics} is dedicated to numerical results, again in the one-dimensional case. We show that the Newton algorithm can be efficiently used in the current least-square optimization setting. In particular, in practice, the exact homogenized coefficients cannot be computed, and only a random approximation of them is available. We monitor here how this randomness propagates to the optimal parameters. The detailed extension of the strategy to the two-dimensional case, along with corresponding numerical tests, will be addressed in a future work~\cite{LMOSxx}.

\section{Discrete homogenization theory}
\label{sec:definition}

For the sake of completeness, we recall first, in Sections~\ref{subsec:homog} and~\ref{sec:truncated}, some elements of homogenization for discrete elliptic equations with random coefficients. We refer to~\cite{K83,K87} for seminal contributions on this topic. For homogenization of elliptic partial differential equations (PDEs), we refer to~\cite{engquist-souganidis} for a general, numerically oriented presentation, to the textbooks~\cite{blp,cd,jikov} and to the review article~\cite{ACLBLT11}. 

Next, in Section~\ref{sec:amael}, we describe the physical background that motivates this work. We eventually turn in Section~\ref{sec:oneD} to the one-dimensional case, where explicit formulas can be obtained.

\subsection{Homogenization result}
\label{subsec:homog}

We first recall some definitions useful for stochastic homogenization, before turning to the specific case of discrete elliptic equations.

Throughout this article, $(\Omega, {\mathcal F}, \xP)$ is a probability space and we denote by $\dps \E(X) = \int_\Omega X(\omega) d\xP(\omega)$ the expectation value of any random variable $X \in \xLone(\Omega, d\xP)$. We next fix $d \in \xN^\star$ (the ambient physical dimension), and assume that the group $(\xZ^d, +)$ acts on $\Omega$. We denote by $(\tau_k)_{k\in \xZ^d}$ this action, and assume that it preserves the measure $\xP$, that is, for all $\displaystyle k \in \xZ^d$ and all $A \in \mathcal{F}$, $\displaystyle \xP(\tau_k A) = \xP(A)$. We assume that the action $\tau$ is {\em ergodic}, that is, if $A \in {\mathcal F}$ is such that $\tau_k A = A$ for any $k \in \xZ^d$, then $\xP(A) = 0$ or 1. 
In addition, we introduce the following notion of stationarity:
\begin{dfntn}
We say that a function $\psi:\xZ^d \times \Omega \to \xR$ is \emph{stationary} if
\begin{equation}
\label{eq:def_stat} 
\forall x, z \in \xZ^d, \quad 
\psi(x+z,\omega)=\psi(x,\tau_z \omega)
\quad \text{a.s.}
\end{equation}
\end{dfntn}

We now focus on the case of discrete elliptic equations. We view $\xZ^d$ as a lattice, whose unit vectors are denoted by $e_i,$ $i \in \{1, \dots, d\}$. Each vertex $x \in \xZ^d$ of the lattice is connected to $2d$ other vertices: $x \pm e_i,$ $i \in \{1, \dots, d\}$. We write $x \sim y$ if $x$ and $y$ are neighbours (i.e. connected), and $e=(x,y)$ the corresponding (non-oriented) edge. For any vertex $x \in \xZ^d$ and any direction $1 \leq i \leq d$, we denote by $a_i(x,\omega) \in (0,\infty)$ the random conductance of the edge $(x,x+e_i)$. We next introduce the diagonal matrix $A$ defined for any vertex $x \in \xZ^d$ by
\begin{equation}
\label{eq:def_A}
A(x,\omega) = \text{diag} \Big( a_1(x,\omega), \ldots, a_d(x,\omega) \Big).
\end{equation}
We assume that, for any direction $i$, the conductances $\left\{ a_i(x,\cdot) \right\}_{x \in \xZ^d}$ form an i.i.d. sequence of random variables. The matrix $A$ is therefore stationary.

\medskip

We next introduce discrete differential operators on the lattice $\xZ^d$. 
\begin{dfntn}
For a function $g : \xZ^d \to \xR$, the gradient $\nabla g : \xZ^d \to \xR^d$ is defined by 
$$
(\nabla g)(x)=\begin{pmatrix} g(x+e_1)-g(x) \\ \vdots \\ g(x+e_d)-g(x) \end{pmatrix}.
$$
For a function $G=(G_1,\dots, G_d): \xZ^d \to \xR^d$, the function $\nabla^\star G : \xZ^d \to \xR$ is defined by 
$$
-(\nabla^\star G)(x)=\sum_{i=1}^d \big( G_i(x)-G_i(x-e_i) \big).
$$
\end{dfntn}
We think of $\nabla^\star G$ as the negative divergence of $G$. The operator $\nabla^\star$ is the $\ell^2$ transpose of $\nabla$ in the following sense: for any compactly supported functions $g:\xZ^d\to\xR$ and $G:\xZ^d\to\xR^d$, 
$$
\sum_{x\in\xZ^d} g(x)\nabla^\star G(x)=\sum_{x \in \xZ^d} \nabla g(x) \cdot G(x).
$$
Hereafter, the notation $a \cdot b$ stands for the usual scalar product in $\xR^d$. 

\medskip

We additionally define rescaled discrete differential operators as follows:
\begin{dfntn}
For a function $g : \eps \xZ^d \to \xR$, the gradient $\nabla_\eps g : \eps \xZ^d \to \xR^d$ is defined by 
$$
(\nabla_\eps g)(x)=\frac{1}{\eps} \begin{pmatrix} g(x+\eps e_1)-g(x) \\ \vdots \\ g(x+\eps e_d)-g(x) \end{pmatrix}.
$$
For a function $G=(G_1,\dots, G_d): \eps \xZ^d \to \xR^d$, the function $\nabla^\star_\eps G : \eps \xZ^d \to \xR$ is defined by 
$$
-(\nabla^\star_\eps G)(x)=\sum_{i=1}^d \frac{G_i(x)-G_i(x-\eps e_i)}{\eps}.
$$
\end{dfntn}

\medskip

The matrix field $A$ is often assumed to satisfy the following assumption:
\begin{hyp}[Ellipticity -- boundedness condition]
\label{ass:ellipticity}
There exist two positive deterministic constants $c$ and $C$ such that the matrix $A$ defined by~\eqref{eq:def_A} satisfies
\begin{equation}
\label{eq:ellipticity}
\forall \xi \in \xR^d, \quad \forall x \in \xZ^d, \quad
c |\xi|^2 \leq \xi \cdot A(x,\omega) \xi \leq C |\xi|^2
\quad \text{a.s.}
\end{equation}
\end{hyp}
In view of~\eqref{eq:def_A}, note that this simply means that $0 < c \leq a_j(x,\omega) \leq C$ almost surely, for any $1 \leq j \leq d$ and any $x \in \xZ^d$.

\medskip

The following homogenization result holds (we refer to~\cite[Theorems~3 and~4]{K83} for a proof):
\begin{thrm}
\label{thm:homogenization}
Let $\D$ be a bounded domain of $\xR^d$ and $f \in \xCzero(\overline{\D})$. Let $A$ be the random stationary matrix field given by~\eqref{eq:def_A}. We assume that~\eqref{eq:ellipticity} holds. Let $u_\eps \in \ell^2(\eps \xZ^d; \xR)$ be the unique solution to
\begin{equation}
\label{eq:pb_base_eps}
\nabla^\star_\eps \big[ A(x/\eps,\omega) \nabla_\eps u_\eps(x,\omega)\big] = f(x) \ \ \text{in $\D \cap \eps \xZ^d$},
\qquad \text{$u_\eps(x,\omega) = 0$ in $(\xR^d \setminus \D) \cap \eps \xZ^d$}.
\end{equation}
When $\eps \to 0$, $u_\eps(\cdot,\omega)$ converges to a homogenized solution $u^\star$ in the following sense.

For any $\xi \in \xR^d$, introduce the \emph{corrector} $\varphi_\xi$ in the direction $\xi$ as the unique solution (defined on $\xZ^d \times \Omega$) to
\begin{equation}
\label{eq:corrector}
\begin{cases} 
-\nabla^\star\ \Big[A(\cdot,\omega) \big(\xi + \nabla \varphi_\xi(\cdot,\omega)\big)\Big]=0 \text{ in $\xZ^d$, a.s.},
\\
\nabla \varphi_\xi \text{ is stationary in the sense of~\eqref{eq:def_stat},} 
\\ 
\forall x \in \xZ^d, \quad \E[\nabla \varphi_\xi(x,\cdot)]=0, 
\\ 
\varphi_\xi(0, \omega)=0  \text{ a.s.}
\end{cases}
\end{equation} 
Introduce next the constant matrix $A^\star$ defined by
\begin{equation}
\label{eq:ahom}
\forall \xi \in \xR^d, \quad 
A^\star\xi = \E\big[ A(x,\cdot) (\xi+\nabla \varphi_\xi(x,\cdot))\big]
\end{equation} 
and the unique solution $u^\star \in \xHone_0(\D)$ to the (continuous) PDE
$$
-{\rm div} \big[ A^\star \widehat{\nabla} u^\star \big] = f \qquad \text{in $\D$},
$$
where $\widehat{\nabla}$ and ${\rm div}$ are the usual (continuous) gradient and divergence differential operators.

Then, we have the (strong) convergence $u_\eps \xrightarrow[\eps \to 0]{} u^\star$, in the sense that 
\begin{equation}
\label{eq:hom_limit}
\eps^d \sum_{x \in \D \cap \eps \xZ^d} | u_\eps(x, \omega) - u^\star(x) |^2 \xrightarrow[\eps \to 0]{} 0 \quad \text{almost surely}.
\end{equation}
\end{thrm}
Note that, in the right-hand side of~\eqref{eq:ahom}, the vector $A(\xi+\nabla \varphi_\xi)$ is stationary, and therefore the expectation may be evaluated at any $x \in \xZ^d$. Note also that, in general, $\varphi_\xi$ itself is not stationary, as the one-dimensional case shows. Only its gradient is. 

\begin{rmrk}
We can define, on $\D$, the function
$$
\widetilde{u}_\eps(x, \omega) = \sum_{k \in \eps \xZ^d \cap \D} u_\eps(k, \omega) \ \mathds{1}_{k + \eps Q}(x), \qquad \text{where $Q=(0,1)^d$}.
$$ 
Then $\widetilde{u}_\eps(\cdot, \omega) \xrightarrow[\eps \to 0]{} u^\star$ in $\xLtwo(\D)$ almost surely.
\end{rmrk}

\subsection{Approximation on finite boxes}
\label{sec:truncated}

The corrector problem~\eqref{eq:corrector} is untractable in practice, since it is posed in the entire lattice $\xZ^d$. Approximations are therefore in order. The standard procedure amounts to considering finite boxes (see e.g.~\cite{BP04}). For a positive integer $N$, we denote by $\Lambda_N$ the finite box $\{0,\dots,N\}^d$ and by $\mathcal{E}_N$ the set of edges in $\Lambda_N$ (see Figure~\ref{fig:finite_box}).

\begin{figure}[H]
  \centering
  \def\Lnet{1.3}    
  \begin{tikzpicture}                               
    \draw[style=help lines] (-1.5*\Lnet,-1.5*\Lnet) grid[step=\Lnet] (4.5*\Lnet,4.5*\Lnet);
    %  Draws the dashed grid
    \foreach \x in {-1,0,...,4} % Two indices running over each
    {
      \foreach \y in {-1,0,...,4} % node on the grid we have drawn 
      {
        \node[draw,circle,inner sep=2pt,fill] at (\Lnet*\x,\Lnet*\y) {};  % Places a dot at those points
        %\draw (\Lnet*\x-0.257,\Lnet*\y-0.5) node[anchor=south]{\pig{0.7}}; % Places a pig at those points
      }
    }
   \draw [thick,-latex,red] (-\Lnet,-\Lnet) -- (-\Lnet,0) node [above left] {$e_2$}; %%%%%% e2
    \draw [thick,-latex,red] (-\Lnet,-\Lnet) -- (0,-\Lnet) node [below right] {$e_1$}; %%%%%% e1

    \filldraw[fill=gray, fill opacity=0.2, draw=black] (-\Lnet*0.5,-\Lnet*0.5)
        rectangle (-\Lnet*0.5+4*\Lnet,-\Lnet*0.5+4*\Lnet);  % gray rectangle --> the finite box
    
    \draw(1*\Lnet,2*\Lnet)node[above right]{$x$}; %%%%% x
    \draw(2*\Lnet,2*\Lnet)node[above right]{$y$}; %%%%% y
    \draw[thick, blue](1*\Lnet,2*\Lnet) -- (2*\Lnet,2*\Lnet);
    \draw[blue](1.5*\Lnet,2*\Lnet) node[below]{$\;e$};   %%%%% e
    \draw(1.5*\Lnet,-0.45*\Lnet)node[below]{$\Lambda_{N}$};   %%%%% N+1

\end{tikzpicture}
  \caption{Finite box $\Lambda_N$ in $\xZ^2$}
  \label{fig:finite_box}
\end{figure}

The \emph{truncated corrector} $\varphi_\xi^N$ defined on $\Lambda_N \times \Omega$ is the unique solution to 
\begin{equation}
\label{eq:trunc_corrector}
\begin{cases}
-\nabla^\star\Big[A(\cdot,\omega)\big( \xi + \nabla\varphi_\xi^N(\cdot,\omega)\big)\Big]=0 \text{ in $\Lambda_N$, a.s.}
\\
\varphi_\xi^N(\cdot,\omega) \text{ is $\Lambda_N$-periodic},
\\ 
\varphi_\xi^N(0, \omega)=0  \text{ a.s.}
\end{cases}
\end{equation}
The homogenized matrix $A^\star$, which is deterministic, is then approximated by the matrix $A^\star_N$ defined by
\begin{equation}
\label{eq:approx_hom}
\forall \xi \in \xR^d, \quad
A^\star_N(\omega) \xi = \frac{1}{| \Lambda_N |} \sum_{x \in \Lambda_N} A(x,\omega)\big(\xi + \nabla \varphi_\xi^N(x,\omega)\big).
\end{equation}
Because of truncation, the practical approximation $A^\star_N$ is random. In the large $N$ limit, the deterministic value is attained, thanks to ergodicity. More precisely, $A^\star_N(\omega)$ converges almost surely towards $ A^\star$ as $N$ goes to infinity, thanks to the ergodic theorem.

\begin{rmrk}
In~\eqref{eq:trunc_corrector}, we have complemented the elliptic equation in $\Lambda_N$ with {\em periodic} boundary conditions. Other choices could be made, such as imposing homogeneous Dirichlet boundary conditions: $\varphi_\xi^N(\cdot,\omega)=0$ on $\partial \Lambda_N$ (see e.g.~\cite{BP04} for a similar discussion in the case of continuous PDEs). In the numerical experiments of Section~\ref{sec:numerics}, we only use periodic boundary conditions, following~\eqref{eq:trunc_corrector}. 
\end{rmrk}

In practice, we work on a finite box $\Lambda_N$, on which the apparent homogenized matrix $A^\star_N$ is random. It is therefore natural to introduce $M$ i.i.d. realizations of the random field $A(x,\omega)$ and solve~\eqref{eq:trunc_corrector}--\eqref{eq:approx_hom} for each of them, thereby obtaining i.i.d. realizations $A^{\star,m}_N(\omega)$, $1 \leq m \leq M$. We next introduce the empirical mean
\begin{equation}
\label{eq:empirical_mean}
\overline{A}^\star_{N,M}(\omega) = \frac{1}{M} \sum_{m=1}^M A^{\star,m}_N(\omega)
\end{equation}
which is, according to the Central Limit Theorem, a converging approximation of $\E \left[ A^\star_N \right]$. We have that
$$
\overline{A}^\star_{N,M}(\omega) \xrightarrow[M\to\infty]{} \E\left[A^\star_N\right] \quad \text{a.s.}
$$
In addition, for any entry $1 \leq i,j \leq d$ of the matrix, we have that, with a probability of 95 \%,
$$
\left| \left( \overline{A}^\star_{N,M}(\omega) \right)_{ij} - \E \left[ \left( A^\star_N \right)_{ij} \right] \right|
\leq
1.96 \sqrt{\frac{\Var \, \left( A^\star_N \right)_{ij}}{M}}.
$$
The error when approximating $A^\star$ by $\overline{A}_{N,M}^\star$ can be written as the sum of two contributions, 
\begin{equation}
\label{eq:errors}
A^\star - \overline{A}_{N,M}^\star= \Big(A^\star-\E[A_N^\star]\Big)+\Big(\E[A_N^\star]-\overline{A}_{N,M}^\star\Big).
\end{equation}
The second term in the right-hand side of~\eqref{eq:errors} is the \emph{statistical} error. The first term is the \emph{systematic} error, due to the fact that, for any finite $N$, $\E[A_N^\star] \neq A^\star$. The dominated convergence theorem ensures that this error vanishes as $N \to \infty$. Many studies have been recently devoted to proving sharp estimates on the rate of this convergence, following the seminal work~\cite{BP04}. In~\cite[Lemma 2.3]{GNO13}, the authors show that the systematic error is of order $N^{-1}$ when the corrector problem is complemented with homogeneous Dirichlet boundary conditions on $\partial \Lambda_N$, and of order $N^{-d} \ln^d(N)$ when using periodic boundary conditions (namely, solving~\eqref{eq:trunc_corrector}). 

\begin{rmrk}
The estimator $\overline{A}^\star_{N,M}$ only agrees with $\E \left[ A^\star_N \right]$ in the limit of an asymptotically large number $M$ of realizations. Note that variance reduction approaches have been introduced in this context (see e.g.~\cite{banff,mprf,cedya} and also~\cite{LM13} for the extension to a nonlinear setting) to obtain approximations of $\E \left[ A^\star_N \right]$ in a more efficient manner than by using $\overline{A}^\star_{N,M}$. 
\end{rmrk}

In the sequel, we will identify the parameters of the microscopic probability distribution on the basis of two types of macroscopic quantities:
\begin{enumerate}
\item the homogenized permeability, which is in practice approximated by $\overline{A}^\star_{N,M}$;
\item the relative variance of any entry $\left( A^\star_N(\omega) \right)_{ij}$, defined by
$$
\VarR \left[ \left( A^\star_N \right)_{ij} \right] 
:= 
\frac{\Var \left[ \left( A^\star_N \right)_{ij} \right]}{\left( \E \left[ \left( A^\star_N \right)_{ij} \right] \right)^2},
$$ 
which is in practice approximated by
\begin{equation}
\label{eq:empirical_var}
S_{N,M} =\frac{1}{(\overline{A}^\star_{N,M})_{ij}^2} \left(\frac{1}{M}\sum_{m=1}^M \left( \left( A^{\star,m}_N(\omega) \right)_{ij} - \left( \overline{A}^\star_{N,M} \right)_{ij} \right)^2\right).
\end{equation}
\end{enumerate}

\subsection{Physical problem}
\label{sec:amael}

We describe here the physical background which inspires this work. As pointed out above, from a physical viewpoint, understanding the microscopic properties of charged porous media is of great importance. Such materials have elaborate geometries that make direct computations very challenging. To circumvent this issue, we use here the Pore Network Model (PNM), which involves a simplified model of the geometry. In the PNM model, pores are located at the vertices of the lattice $\xZ^d$. Neighbouring pores are connected by channels, which allow water to flow. Each channel $(x,x+e_i)$ is endowed with its random conductance $a_i(x,\omega)>0$, the probability distribution of which is discussed below.

Experiments provide measures on the \emph{macroscopic permeability} $K_{\rm obs}$, which is modelled as a homogenized coefficient $K^\star$. In practice, as explained in Section~\ref{sec:truncated}, the homogenized coefficient can only be approximated through a computation on a large box. Assuming that the conductance field $a_i(x,\omega)$ is given for any direction $1 \leq i \leq d$ and any vertex $x$ on the finite lattice $\Lambda_N$, the PNM model consists in computing the pressure field $P(x,\omega)$ by solving the conservation equations (i.e., Darcy law) in the network. This leads to the following linear system: 
\begin{equation}
\label{eq:conservation}
\forall x \in \Lambda_N, \quad 
\sum_{y \sim x} \tilde{a}(x,y,\omega) \big(P(y,\omega)-P(x,\omega)\big)
=0,
\end{equation}
where $\tilde{a}(x,y,\omega)$ is the conductance of the non-oriented edge $(x,y)$. Some boundary conditions need to be imposed to make this problem well-posed, they are discussed below. We next see, by definition of $a_i$, that
\begin{eqnarray}
&&
\sum_{y \sim x} \tilde{a}(x,y,\omega) \big(P(y,\omega)-P(x,\omega)\big)
\nonumber
\\
&=&
\sum_{i=1}^d a_i(x,\omega) \big(P(x+e_i,\omega)-P(x,\omega)\big)
+
\sum_{i=1}^d a_i(x-e_i,\omega) \big(P(x-e_i,\omega)-P(x,\omega)\big)
\nonumber
\\
&=&
\nabla^\star \big[ A(\cdot,\omega)\nabla P(\cdot,\omega) \big](x),
\label{eq:recast}
\end{eqnarray}
where the matrix $A$ is defined in terms of $\left\{ a_i \right\}_{i=1}^d$ by~\eqref{eq:def_A}. 

\medskip

\begin{figure}[H]
  \centering
  \def\Lnet{1.3}    
  \begin{tikzpicture}
                                   
    \draw[style=help lines] (0,0) grid[step=\Lnet] (3*\Lnet,3*\Lnet);
    %  Draws the dashed grid   
    \foreach \x in {0,...,3} % Two indices running over each
    {
      \foreach \y in {0,...,3} % node on the grid we have drawn 
      {
        \node[draw,circle,inner sep=2pt,fill] at (\Lnet*\x,\Lnet*\y) {}; % Places a dot at those points
      }
    }  
    \node[draw, circle, inner sep=2pt,fill] at (-\Lnet,\Lnet*1.5) {}  ;
    
    \draw (-\Lnet,\Lnet*1.5) node [above left] {$P_O$};
    
    \foreach \y in {0,...,3}
    {
      \draw [dotted] (-\Lnet,\Lnet*1.5) -- (0,\Lnet*\y) {};
    }
    
    \node[draw, circle, inner sep=2pt,fill] at (4*\Lnet,1.5*\Lnet) {};    
    \draw  (4*\Lnet,1.5*\Lnet) node [above right] {$P_I$};
    \foreach \y in {0,...,3}
    {
      \draw [dotted] (4*\Lnet,\Lnet*1.5) -- (3*\Lnet,\Lnet*\y) {};
    }    
    \draw [thick,-latex,olive] (0.65*\Lnet,0.5*\Lnet) -- (1.65*\Lnet,0.5*\Lnet) node [above left] {$\xi\;\;$}; %%%%%% e2
    
    \draw [thick,-latex,red] (-\Lnet,-\Lnet) -- (-\Lnet,0) node [above left] {$e_2$}; %%%%%% e2
    \draw [thick,-latex,red] (-\Lnet,-\Lnet) -- (0,-\Lnet) node [below right] {$e_1$}; %%%%%% e1
    \draw(1*\Lnet,2*\Lnet)node[above right]{$x$}; %%%%% x
    \draw(2*\Lnet,2*\Lnet)node[above right]{$y$}; %%%%% y
    \draw[thick, blue](1*\Lnet,2*\Lnet) -- (2*\Lnet,2*\Lnet);
    \draw[blue](1.5*\Lnet,2*\Lnet) node[below]{$\;e$};   %%%%% e
    \draw(1.5*\Lnet,-0.45*\Lnet)node[below]{$\Lambda_{N}$};   %%%%% N+1
  \end{tikzpicture}
  \caption{Finite lattice with boundary conditions}
  \label{fig:finite_box_bc}
\end{figure}
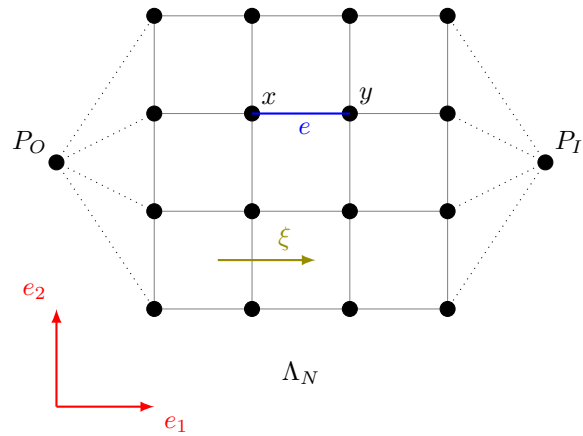

We now describe (in the two-dimensional case, for the sake of simplicity) the boundary conditions imposed on~\eqref{eq:conservation}. They are designed to mimic experimental conditions. We first recall that the large box reads $\Lambda_N = \{0,\dots,N\}^2$. The pressure field is assumed to be periodic in the vertical direction, whereas a macroscopic gradient is imposed in the horizontal direction as follows. Imagine that all vertices with coordinates $(0,\cdot)$ are connected to one fixed vertex denoted by $O$, representing a pressure reservoir at pressure $P_O$. Likewise, all vertices with coordinates $(N,\cdot)$ are connected to one fixed vertex denoted by $I$ at pressure $P_I$ (see Figure~\ref{fig:finite_box_bc}). Then, the boundary conditions write 
$$
\text{ for all $j\in\{0,\dots,N\}$}, \quad
P(0,j)=P_O \ \ \text{and} \ \ P(N,j)=P_I.
$$
Once~\eqref{eq:conservation} is solved with the above boundary conditions, the \emph{macroscopic permeability} $K^\star_N$ is defined by
\begin{equation}
\label{eq:def_k_star}
K^\star_N(\omega) := \frac{N}{P_O-P_I} \frac{1}{|\Lambda_N|} \sum_{x \in \Lambda_N} \tilde{a}(x,x+e_1,\omega) \big(P(x,\omega)-P(x+e_1,\omega)\big). 
\end{equation}

Let us now show that Equations~\eqref{eq:conservation}--\eqref{eq:def_k_star} are actually the same as Equations~\eqref{eq:trunc_corrector}--\eqref{eq:approx_hom} written above. 

By linearity of~\eqref{eq:conservation}--\eqref{eq:def_k_star}, we can always assume that $P_O = 0$ and $P_I = N$. Let $P$ be a solution to~\eqref{eq:conservation}. We introduce $\varphi_{e_1}$ such that
$$
P(x,\omega) = x \cdot e_1 +\varphi_{e_1}(x,\omega).
$$
In view of~\eqref{eq:conservation} and~\eqref{eq:recast}, we see that $\varphi_{e_1}(\cdot,\omega)$ is solution to
$$
\forall x \in \Lambda_N, \quad 
\nabla^\star \big[ A(\cdot,\omega) (e_1 + \nabla \varphi_{e_1}(\cdot,\omega) \big](x) = 0
$$
with $\varphi_{e_1}((0,j),\omega) = \varphi_{e_1}((N,j),\omega) = 0$ for any $j$ and $\varphi_{e_1}(\cdot,\omega)$ is periodic in the vertical direction. Up to the choice of boundary conditions, we thus recognize~\eqref{eq:trunc_corrector} for $\xi = e_1$. We also infer from~\eqref{eq:def_k_star} that
\begin{eqnarray*}   
K^\star_N(\omega) 
&=& 
\frac{N}{P_O-P_I} \frac{1}{|\Lambda_N|} \sum_{x \in \Lambda_N} \tilde{a}(x,x+e_1,\omega) \big(P(x,\omega)-P(x+e_1,\omega)\big)
\\
&=&
\frac{1}{|\Lambda_N|} \sum_{x \in \Lambda_N} a_1(x,\omega) \big(\varphi_{e_1}(x,\omega)-\varphi_{e_1}(x+e_1,\omega) - e_1 \cdot e_1 \big)
\\
&=&
\frac{1}{|\Lambda_N|} \sum_{x \in \Lambda_N} a_1(x,\omega) \ e_1^T \ \big(e_1 + \nabla \varphi_{e_1}(x,\omega) \big)
\\
&=&
e_1^T A^\star_N(\omega) e_1,
\end{eqnarray*}
where $A^\star_N(\omega)$ is defined by~\eqref{eq:approx_hom}, and where we have used~\eqref{eq:def_A} in the last line. Thus, up to the choice of boundary conditions in the corrector problem, the formulation~\eqref{eq:conservation}--\eqref{eq:def_k_star} is identical to the formulation~\eqref{eq:trunc_corrector}--\eqref{eq:approx_hom}. 

\medskip

\begin{figure}[htbp]
\centering
\includegraphics[scale=0.4]{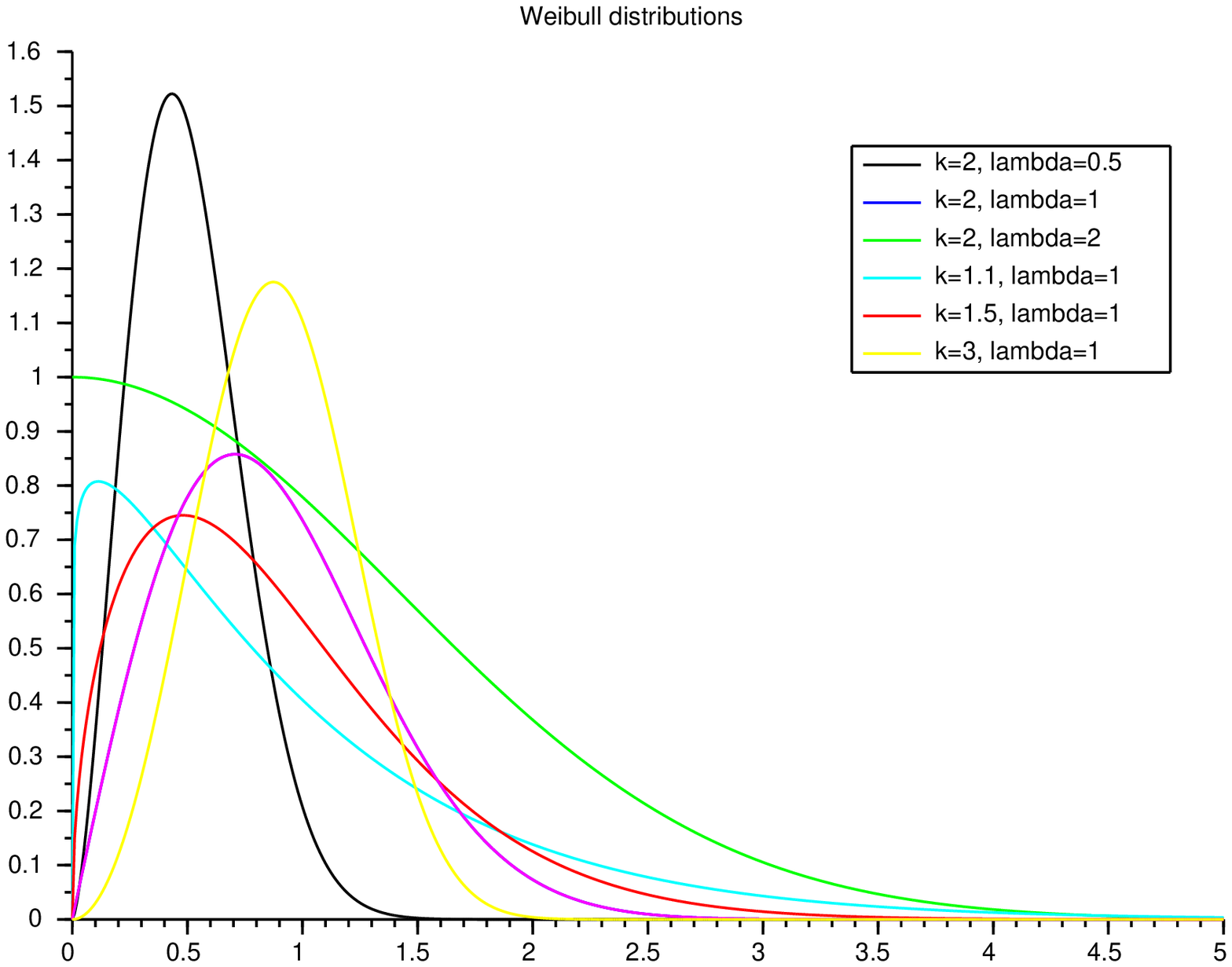}
\caption{Examples of Weibull distributions.}
\label{fig:weibull}
\end{figure}

We eventually discuss the choice of the probability distribution for the conductances. Based on experimental results, it is reasonable to assume the following: 

\begin{hyp}
\label{ass:weibull4} 
We assume that the radius $r$ of the channels are i.i.d. random variables distributed according to a Weibull law of parameter $\theta := (\lambda,k) \in (\xR_+^\star)^2$, that we denote ${\cal W}(\lambda,k)$. We recall that such random variables are positive, with a probability density that reads (see Figure~\ref{fig:weibull})
$$
\forall r > 0, \qquad f(r; k, \lambda) = \frac{k}{\lambda} \left ( \frac{r}{\lambda} \right)^{k-1} \exp \left( -(r/\lambda)^k \right),
$$   
corresponding to the cumulative distribution function
$$
F(r; k, \lambda) = \int_0^r f(s; k, \lambda) \, ds = 1 - \exp \left( -(r/\lambda)^k \right).
$$
Note that the radius of all channels (independently of their direction $1 \leq i \leq d$) share the same probability distribution.  
\end{hyp}
In practice, a Weibull distribution is generated as follows. Let $u(\omega)$ be a random variable uniformly distributed in $[0,1]$. Then
$$
r(\omega) = \lambda \Big[ - \ln (1 - u(\omega)) \Big]^{1/k}
$$
is distributed according to the Weibull law of parameter $(\lambda,k)$.

\medskip

Physical arguments lead to the fact that the conductance $a_i(x,\omega)$ of any channel $(x,x+e_i)$ is directly related to its radius $r(x,x+e_i,\omega)$. Hereafter, we assume that
\begin{equation}
\label{eq:cov}
a_i(x,\omega) = C_0 \, r^4(x,x+e_i,\omega) = C_0 \, \lambda^4 \Big[ - \ln (1 - u(\omega)) \Big]^{4/k},
\end{equation}
where $C_0$ is a constant (for instance, for a Poiseuille flow, $C_0=\pi /(8\eta)$ where $\eta$ is the fluid viscosity). For the sake of simplicity, we will take $C_0 = 1$ in the sequel. Therefore, we assume that
\begin{equation}
\label{eq:prop}
\begin{array}{c}
\text{\emph{The conductances $\left\{ a_i(x,\omega) \right\}_{x \in \xZ^d, \, 1 \leq i \leq d}$ form an i.i.d. sequence of random variables}}
\\ 
\text{\emph{that are distributed according to the Weibull law of parameter $(\lambda^4,k/4)$.}}
\end{array}
\end{equation}

\begin{rmrk}
Note that the Weibull distribution is isolated neither from 0 nor from $\infty$. The above model therefore does not satisfy the ellipticity condition~\eqref{eq:ellipticity}. First, we show in Section~\ref{sec:1D_pas_elliptique} below that, in the one-dimensional case, the assumption~\eqref{eq:ellipticity} is not necessary, and that homogenization holds under a weaker assumption. Second, we refer to~\cite{B11} for similar studies (again under assumptions weaker than~\eqref{eq:ellipticity}) in higher-dimensional cases.
\end{rmrk}

\begin{rmrk}
The numerical tests of Section~\ref{sec:numerics} are performed with the above model, and thus aim at identifying the two parameters $\lambda$ and $k$. We however note that nothing in our approach is specific to this particular model using Weibull laws. This choice is only motivated by physical reasons.
\end{rmrk}

Since the conductances $a_i(x,\omega)$ are all i.i.d. (for any $1 \leq i \leq d$ and any $x \in \xZ^d$), the problem is invariant by any rotation of angle $\pi/2$. The homogenized matrix $A^\star$ is therefore proportional to the identity matrix $\mbox{Id}_d$, and reads
$$
A^\star = K^\star \, \mbox{Id}_d
$$
where $K^\star \in (0,\infty)$ is the homogenized permeability. We can also write that 
$$
K^\star = e_1 \cdot A^\star e_1.
$$
In practice, we only have access to $A^\star_N(\omega)$, which is a symmetric matrix (but is a priori not proportionnal to the identity matrix). Because all directions are statistically identical, we only focus on
\begin{equation}
\label{eq:def_kstar_N}
K^\star_N(\omega) := e_1 \cdot A^\star_N(\omega) e_1.
\end{equation}

\subsection{The one dimensional case}
\label{sec:oneD}

The purpose of this section is two-fold. First, we provide explicit formulas for the homogenized quantities in terms of the microscopic field $A(x,\omega)$. We derive these formulas assuming that~\eqref{eq:ellipticity} holds. Second, we show that we can relax Assumption~\eqref{eq:ellipticity} and still state a homogenization result. 

\subsubsection{Explicit formulas in the elliptic case~\eqref{eq:ellipticity}}

In the one-dimensional case, the problem~\eqref{eq:trunc_corrector}--\eqref{eq:approx_hom} can be analytically solved. We have
\begin{equation}
\label{eq:approx_hom_1D}
A^\star_N(\omega) = \left( \frac 1N \sum_{x \in \Lambda_N} \frac 1 {A(x, \omega)} \right)^{-1} \quad \text{ for almost all $\omega$}.
\end{equation}
Likewise, the problem~\eqref{eq:corrector}--\eqref{eq:ahom} can also be solved, yielding the formula
\begin{equation}
\label{eq:ahom_1D}
A^\star=\left(\E\left[\frac{1}{A(x, \cdot)}\right]\right)^{-1},
\end{equation}
which can be evaluated at any $x \in \xZ$ due to the stationarity of $A$. 

First, it can be checked that the homogenization convergence~\eqref{eq:hom_limit} holds, and second, that $A^\star_N(\omega)$ indeed converges to $A^\star$ when $N \to \infty$.

We note that, as soon as $A(x,\omega) > 0$ a.s. for any $x \in \xZ$ and $A^{-1}(x,\cdot) \in \xLone(\Omega)$ (this latter condition being independent of $x$), formulas~\eqref{eq:approx_hom_1D} and~\eqref{eq:ahom_1D} are well-defined. The aim of the next section is to recall that, in the one-dimensional case, these assumptions are enough for homogenization to hold.

\subsubsection{Relaxing Assumption~\eqref{eq:ellipticity}}
\label{sec:1D_pas_elliptique}

In this section, we show that the following assumption is enough for homogenization to hold:
\begin{hyp}
\label{ass:ellipticity_weak}
We assume that the coefficient $A$ is almost surely positive and finite and satisfies 
\begin{equation}
\label{eq:ellipticity_weak}
A^{-1}(x,\cdot) \in \xLone(\Omega).
\end{equation}
\end{hyp}
Of course, by stationarity, if~\eqref{eq:ellipticity_weak} is satisfied for some $x \in \xZ$, then it is satisfied for all $x \in \xZ$.

\begin{thrm}
\label{thm:homogenization_1D}
Let $\D$ be a bounded domain of $\xR$, $f \in \xCzero(\overline{\D})$ and $A$ be a random stationary scalar field (defined on $\xZ \times \Omega$) that satisfies~\eqref{eq:ellipticity_weak}. Let $u_\eps \in \ell^2(\eps \xZ; \xR)$ be the unique solution to
\begin{equation}
\label{eq:pb_base_eps_1D}
\nabla^\star_\eps \big[ A(x/\eps,\omega) \nabla_\eps u_\eps(x,\omega)\big] = f(x) 
\ \ \text{in $\D \cap \eps \xZ$},
\qquad \text{$u_\eps(x,\omega) = 0$ in $(\xR \setminus \D) \cap \eps \xZ$},
\end{equation}
and let $u^\star \in \xHone_0(\D)$ be the unique solution to the (continuous) boundary value problem
\begin{equation}
\label{eq:pb_homog_1D}
- \big[ A^\star (u^\star)' \big]' = f \qquad \text{in $\D$},
\end{equation}
where $A^\star$ is defined by~\eqref{eq:ahom_1D}.

Then, when $\eps \to 0$, $u_\eps(\cdot,\omega)$ converges to the homogenized solution $u^\star$, in the sense that
\begin{equation}
\label{eq:hom_limit_1D}
\eps \sum_{x \in \D \cap \eps \xZ} | u_\eps(x, \omega) - u^\star(x) |^2 \xrightarrow[\eps \to 0]{} 0 \quad \text{almost surely}.
\end{equation}
\end{thrm}

Note that~\eqref{eq:pb_base_eps_1D} is almost surely well-posed. Indeed, since $A$ is stationary and $0 < A(0,\omega) < \infty$ almost surely, we have that, almost surely, $0 < A(x,\omega) < \infty$ for any $x \in \xZ$. For those $\omega$, problem~\eqref{eq:pb_base_eps_1D} is well-posed. Likewise, since $A$ is almost surely finite (resp. $A^{-1}(0,\cdot) \in \xLone(\Omega)$), we have that $A^\star < \infty$ (resp. $A^\star > 0$) and hence~\eqref{eq:pb_homog_1D} is well-posed.

\begin{proof}
The proof proceeds by truncation of the coefficient $A$ in the neighbourhood of 0 and $+\infty$. For the sake of simplicity, we take $\D = (0,1)$. For any $m \in \xN^\star$, we introduce the coefficient $A_m$ defined on $\xZ \times \Omega$ by
$$
A_m(x,\omega) :=
\left\{
\begin{array}{ccc}
\dps \frac{1}{m} & \quad \text{if} \quad & \dps 0 < A(x,\omega) < \frac{1}{m},
\\
A(x,\omega) & \quad \text{if} \quad & \dps \frac{1}{m} \leq A(x,\omega) \leq m,
\\
m & \quad \text{if} \quad & A(x,\omega) > m.
\end{array}
\right.
$$
We set
$$
%\begin{equation}
%\label{eq:ahom_1D_m}
A^\star_m =\left(\E\left[\frac{1}{A_m(0, \cdot)}\right]\right)^{-1}.
%\end{equation}
$$
For almost all $\omega$ (i.e. those such that $A(0,\omega) > 0$), we have
\begin{gather*}
\lim_{m \to \infty} \frac{1}{A_m(0,\omega)} = \frac{1}{A(0,\omega)},
\\
\forall m \in \xN^\star, \quad 0 < \frac{1}{A_m(0,\omega)} \leq 1 + \frac{1}{A(0,\omega)}, 
\end{gather*}
where the right-hand side of the above second line belongs to $\xLone(\Omega)$, in view of the assumption~\eqref{eq:ellipticity_weak}. Therefore, the dominated convergence theorem implies that
\begin{equation}
\label{eq:william1} 
\lim_{m \to \infty} A^\star_m = A^\star.
\end{equation} 
Let $u^m_\eps \in \ell^2(\eps \xZ; \xR)$ be the unique solution to
\begin{equation}
\label{eq:pb_base_eps_1D_m}
\nabla^\star_\eps \big[ A_m(x/\eps,\omega) \nabla_\eps u^m_\eps(x,\omega)\big] = f(x) 
\ \ \text{in $(0,1) \cap \eps \xZ$},
\qquad \text{$u^m_\eps(x,\omega) = 0$ in $(\xR \setminus (0,1)) \cap \eps \xZ$},
\end{equation}
and let $u_m^\star \in \xHone_0(0,1)$ be the unique solution to the (continuous) boundary value problem
$$
- \big[ A_m^\star (u_m^\star)' \big]' = f \qquad \text{in $(0,1)$}.
$$
We write
\begin{equation}
\label{eq:william2}
\| u_\eps(\cdot,\omega) - u^\star \|_{\ell^2_\eps} 
\leq
\| u_\eps(\cdot,\omega) - u^m_\eps(\cdot,\omega) \|_{\ell^2_\eps} 
+
\| u^m_\eps(\cdot,\omega) - u_m^\star \|_{\ell^2_\eps}
+
\| u_m^\star - u^\star \|_{\ell^2_\eps}
\end{equation}
where, for any function $v$,
$$
\| v \|_{\ell^2_\eps} := \sqrt{ \eps \sum_{x \in (0,1) \cap \eps \xZ} v^2(x) }.
$$
We successively study the three terms of the right-hand side of~\eqref{eq:william2}.

First, we have
$$
\lim_{\eps \to 0} \| u_m^\star - u^\star \|_{\ell^2_\eps} = \| u_m^\star - u^\star \|_{\xLtwo(0,1)},
$$ 
and the convergence~\eqref{eq:william1} implies that
\begin{equation}
\label{eq:william3}
\lim_{m \to \infty} \lim_{\eps \to 0} \| u_m^\star - u^\star \|_{\ell^2_\eps} = 0.
\end{equation}
Second, the coefficient $A_m$ satisfies the ellipticity condition~\eqref{eq:ellipticity}, so we infer from Theorem~\ref{thm:homogenization} that, for any $m \in \xN^\star$,
\begin{equation}
\label{eq:william4}
\lim_{\eps \to 0} \| u^m_\eps(\cdot,\omega) - u_m^\star \|_{\ell^2_\eps} = 0 \quad \text{a.s.}
\end{equation}
We eventually turn to the first term of the right-hand side of~\eqref{eq:william2}. Let
$$
F_\eps(x) = \eps \sum_{y \in (0,x] \cap \eps \xZ} f(y),
$$
which satisfies, for any $x$, $\left| F_\eps(x) \right| \leq \| f \|_{\xLinfty}$.
Integrating once the equations~\eqref{eq:pb_base_eps_1D} and~\eqref{eq:pb_base_eps_1D_m}, we can show that there exist two random variables $C_\eps(\omega)$ and $C^m_\eps(\omega)$, independent of $x$, such that
\begin{eqnarray}
\label{eq:wx1}
A_m \left(\frac{x}{\eps}, \omega\right) \nabla_\eps u^m_\eps (x, \omega) 
&=& 
- F_\eps(x) + C^m_\eps(\omega),
\\
\label{eq:wx2}
A \left(\frac{x}{\eps}, \omega\right) \nabla_\eps u_\eps (x, \omega) 
&=& 
- F_\eps(x) + C_\eps(\omega).
\end{eqnarray}
Using the boundary conditions on $u^m_\eps$ and $u_\eps$, we get
$$
C_\eps(\omega) = \frac{ {\cal N}_\eps(\omega) }{ {\cal D}_\eps(\omega) }
\quad \text{and} \quad
C^m_\eps(\omega) = \frac{ {\cal N}^m_\eps(\omega) }{ {\cal D}^m_\eps(\omega) }
$$
where
\begin{eqnarray*}
{\cal D}_\eps(\omega)
&=&
\eps \sum_{x \in (0,1) \cap \eps \xZ} A\left(\frac{x}{\eps}, \omega \right)^{-1}, 
\qquad
{\cal N}_\eps(\omega)
=
\eps \sum_{x \in (0,1) \cap \eps \xZ} A\left(\frac{x}{\eps}, \omega \right)^{-1} F_\eps(x), 
\\
{\cal D}^m_\eps(\omega)
&=&
\eps \sum_{x \in (0,1) \cap \eps \xZ} A_m\left(\frac{x}{\eps}, \omega \right)^{-1},
\qquad
{\cal N}^m_\eps(\omega)
=
\eps \sum_{x \in (0,1) \cap \eps \xZ} A_m\left(\frac{x}{\eps}, \omega \right)^{-1} F_\eps(x).
\end{eqnarray*}
All these quantities are well-defined for almost all $\omega$. 
We claim that
\begin{equation}
\label{eq:william6}
\lim_{m \to \infty} \limsup_{\eps \to 0}
\left| C_\eps(\omega) - C^m_\eps(\omega) \right|
= 0 \quad \text{a.s.}
\end{equation}
To prove this claim, we start by writing that
\begin{equation}
\label{eq:william5}
C_\eps(\omega) - C^m_\eps(\omega)
=
\frac{ {\cal N}_\eps(\omega) - {\cal N}^m_\eps(\omega) }{ {\cal D}_\eps(\omega) }
+
\frac{ {\cal N}^m_\eps(\omega) }{ {\cal D}^m_\eps(\omega) {\cal D}_\eps(\omega) }
\left( {\cal D}^m_\eps(\omega) - {\cal D}_\eps(\omega) \right).
\end{equation}
Introduce
$$
b_m\left(x, \omega \right) 
= 
\left| \frac{1}{A_m\left(x, \omega \right)} - \frac{1}{A\left(x, \omega \right)}\right|
\quad \text{and} \quad
{\cal B}^m_\eps(\omega)
=
\eps \sum_{x \in (0,1) \cap \eps \xZ} b_m\left(\frac{x}{\eps}, \omega \right).
$$
For any $m \in \xN^\star$, we get
\begin{eqnarray}
\label{eq:ww1}
\left| {\cal N}^m_\eps(\omega) - {\cal N}_\eps(\omega) \right|
&\leq&
\| f \|_{\xLinfty} \ \eps \sum_{x \in (0,1) \cap \eps \xZ} b_m\left(\frac{x}{\eps}, \omega \right)
=
\| f \|_{\xLinfty} \ {\cal B}^m_\eps(\omega),
\\
\label{eq:ww2}
\left| {\cal D}^m_\eps(\omega) - {\cal D}_\eps(\omega) \right|
&\leq&
\eps \sum_{x \in (0,1) \cap \eps \xZ} b_m\left(\frac{x}{\eps}, \omega \right)
=
{\cal B}^m_\eps(\omega),
\\
\label{eq:ww3}
\left| {\cal N}^m_\eps(\omega) \right|
&\leq&
\| f \|_{\xLinfty} \ {\cal D}^m_\eps(\omega).
\end{eqnarray}
Using the ergodic theorem for the stationary functions $A^{-1}$, $A_m^{-1}$ and $b_m$, we have that, for any $m \in \xN^\star$, almost surely,
\begin{equation}
\label{eq:ww4}
\lim_{\eps \to 0} {\cal D}_\eps(\omega) = \frac{1}{A^\star}, \quad
\quad
\lim_{\eps \to 0} {\cal D}^m_\eps(\omega) = \frac{1}{A_m^\star}, \quad
\quad
\lim_{\eps \to 0} {\cal B}^m_\eps(\omega) = {\cal B}^m_\star := \E\left[ \left| \frac{1}{A_m(0, \cdot)} - \frac{1}{A(0, \cdot)} \right| \right].
\end{equation}
We introduce
$$
\Omega_{\rm conv} = \left\{
\omega \in \Omega; \
\lim_{\eps \to 0} {\cal D}_\eps(\omega) = \frac{1}{A^\star} \
\ \text{and, for all $m \in \xN^\star$}, \
\lim_{\eps \to 0} {\cal D}^m_\eps(\omega) = \frac{1}{A_m^\star}, \
\quad
\lim_{\eps \to 0} {\cal B}^m_\eps(\omega) = {\cal B}^m_\star
\right\}
$$
and we deduce that $\xP(\Omega_{\rm conv}) = 1$.

Let $\omega \in \Omega_{\rm conv}$. In view of~\eqref{eq:ww4}, we know that there exists $\eps^m_0(\omega)$ such that, for any $\eps < \eps^m_0(\omega)$, we have
\begin{equation}
\label{eq:ww5}
\frac{1}{2A^\star} \leq {\cal D}_\eps(\omega),
\quad
\frac{1}{2A_m^\star} \leq {\cal D}^m_\eps(\omega) \leq \frac{3}{2A_m^\star}.
\end{equation}
We thus infer from~\eqref{eq:william5}, \eqref{eq:ww5}, \eqref{eq:ww1}, \eqref{eq:ww3} and~\eqref{eq:ww2} that, for any $\omega \in \Omega_{\rm conv}$, any $m \in \xN^\star$ and any $\eps < \eps^m_0(\omega)$, we have
\begin{eqnarray}
\nonumber
\left| C_\eps(\omega) - C^m_\eps(\omega) \right|
&\leq&
2 A^\star \left| {\cal N}_\eps(\omega) - {\cal N}^m_\eps(\omega) \right| 
+
4 A^\star A^\star_m \left| {\cal N}^m_\eps(\omega) \right| \
\left| {\cal D}^m_\eps(\omega) - {\cal D}_\eps(\omega) \right|
\\
\nonumber
&\leq&
2 A^\star \ \| f \|_{\xLinfty} \ {\cal B}^m_\eps(\omega)
+
4 A^\star A^\star_m \ \| f \|_{\xLinfty} \ {\cal D}^m_\eps(\omega) \
{\cal B}^m_\eps(\omega)
\\
\label{eq:wx3}
&\leq&
2 A^\star \ \| f \|_{\xLinfty} \ {\cal B}^m_\eps(\omega)
+
6 A^\star \ \| f \|_{\xLinfty} \ {\cal B}^m_\eps(\omega).
\end{eqnarray}
Hence, for any $\omega \in \Omega_{\rm conv}$ and any $m \in \xN^\star$, we have
$$
\limsup_{\eps \to 0}
\left| C_\eps(\omega) - C^m_\eps(\omega) \right|
\leq
8 A^\star \ \| f \|_{\xLinfty} \ {\cal B}^\star_m.
$$
The dominated convergence theorem implies that $\dps \lim_{m \to \infty} {\cal B}^\star_m = 0$, hence, for any $\omega \in \Omega_{\rm conv}$, we have
$$
\lim_{m \to \infty} \limsup_{\eps \to 0}
\left| C_\eps(\omega) - C^m_\eps(\omega) \right|
= 0.
$$
Since $\xP(\Omega_{\rm conv}) = 1$, we have proved the claim~\eqref{eq:william6}.

\medskip

We now proceed and deduce from~\eqref{eq:wx1} and~\eqref{eq:wx2} that
\begin{eqnarray*}
u^m_\eps(z, \omega) 
&=& 
\eps \sum_{x \in (0,z) \cap \eps \xZ} A_m \left(\frac{x}{\eps}, \omega\right)^{-1} \left( C^m_\eps(\omega) - F_\eps(x)\right),
\\
u_\eps(z, \omega) 
&=& 
\eps \sum_{x \in (0,z) \cap \eps \xZ} A \left(\frac{x}{\eps}, \omega\right)^{-1} \left( C_\eps(\omega) - F_\eps(x)\right),
\end{eqnarray*}
hence
$$
\left| u^m_\eps(z, \omega) - u_\eps(z, \omega) \right|
\leq
\left| C^m_\eps(\omega) - C_\eps(\omega) \right| \ {\cal D}^m_\eps(\omega)
+
\left( \left| C_\eps(\omega) \right| + \| f \|_{\xLinfty} \right) \ {\cal B}^m_\eps(\omega).
$$
Using that $\left| C_\eps(\omega) \right| \leq \| f \|_{\xLinfty}$, we deduce that 
$$
\| u^m_\eps(\cdot,\omega) - u_\eps(\cdot,\omega) \|_{\ell^2_\eps} 
\leq
\left| C^m_\eps(\omega) - C_\eps(\omega) \right| \ {\cal D}^m_\eps(\omega)
+
2 \| f \|_{\xLinfty} \ {\cal B}^m_\eps(\omega).
$$
For any $\omega \in \Omega_{\rm conv}$, any $m \in \xN^\star$ and any $\eps < \eps^m_0(\omega)$, using~\eqref{eq:wx3} and~\eqref{eq:ww5}, we obtain that
$$
\| u^m_\eps(\cdot,\omega) - u_\eps(\cdot,\omega) \|_{\ell^2_\eps} 
\leq
8 A^\star \ \| f \|_{\xLinfty} \ {\cal B}^m_\eps(\omega) \ \frac{3}{2 A^\star_m}
+
2 \| f \|_{\xLinfty} \ {\cal B}^m_\eps(\omega),
$$
hence, for any $\omega \in \Omega_{\rm conv}$ and any $m \in \xN^\star$,
$$
\limsup_{\eps \to 0} \| u^m_\eps(\cdot,\omega) - u_\eps(\cdot,\omega) \|_{\ell^2_\eps}
\leq
8 A^\star \ \| f \|_{\xLinfty} \ {\cal B}^m_\star \ \frac{3}{2 A^\star_m}
+
2 \| f \|_{\xLinfty} \ {\cal B}^m_\star,
$$
and thus, almost surely, 
\begin{equation}
\label{eq:william7}
\lim_{m \to \infty} \limsup_{\eps \to 0} \| u^m_\eps(\cdot,\omega) - u_\eps(\cdot,\omega) \|_{\ell^2_\eps}
=
0.
\end{equation}
Collecting~\eqref{eq:william2}, \eqref{eq:william3}, \eqref{eq:william4} and~\eqref{eq:william7}, we obtain that
$$
\limsup_{\eps \to 0} 
\| u_\eps(\cdot,\omega) - u^\star \|_{\ell^2_\eps} 
=0 \quad \text{a.s.,}
$$
which is the convergence~\eqref{eq:hom_limit_1D}.
\end{proof}

\subsubsection{The case of Weibull laws}

Following Section~\ref{sec:amael}, assume that the conductances are given by~\eqref{eq:prop}, i.e. are distributed according to the Weibull law of parameter $(\lambda^4,k/4)$. For any $k>0$, Assumption~\eqref{eq:ellipticity} is not satisfied. However, when $k>4$, Assumption~\eqref{eq:ellipticity_weak} is satisfied and we have, in view of Theorem~\ref{thm:homogenization_1D}, that
\begin{equation}
\label{eq:astar_1D_weibull}
A^\star=\frac{\lambda^4}{\Gamma(1-4/k)}
\end{equation}
where $\Gamma$ is the Euler Gamma function defined for any $z>0$ by 
\begin{equation}
\label{eq:def_Gamma}
\Gamma(z) = \int_0^\infty t^{z-1} \, \exp(-t) \, dt.
\end{equation}
The variance of $A^\star_N$ is finite if and only if $k>8$. In the sequel, we work in the range $k > 8$.

\section{A parameter fitting problem} 
\label{sec:inverse}

We now describe the problem we consider, first in the general case (Section~\ref{sec:inv_dD}), next in the one-dimensional case (Section~\ref{sec:inv_1D}). In that latter section, we also motivate our choice of macroscopic quantities from which we fit the parameters of the Weibull law.

\subsection{General case}
\label{sec:inv_dD}

We assume that we are given two observed quantities, the first coefficient in the macroscopic permeability matrix (see~\eqref{eq:def_kstar_N})
$$
K^{\star,{\rm obs}}_N(\omega) = e_1 \cdot A^{\star,{\rm obs}}_N(\omega) e_1
$$
and its relative variance $S^{\rm obs}_N$ for some parameters $\theta_{\rm obs}=(\lambda_{\rm obs},k_{\rm obs})$ of the Weibull law. Note that the relative variance crucially depends on the size $N^d$ of the finite box on which it is measured (in contrast to the apparent permeability, which converges to a finite value when $N \to \infty$). We assume here that we know this size. In practice, these three quantities, $N$, $K^{\star,{\rm obs}}_N$ and $S^{\rm obs}_N$, can be obtained by physical experiments. We therefore assume that there exists $\theta_{\rm obs}$ and $N$ such that
\begin{equation}
\label{eq:synt}
\E[ K^\star_N(\cdot,\theta_{\rm obs}) ] = K^{\star,{\rm obs}}_N, 
\qquad
\VarR[ K^\star_N(\cdot,\theta_{\rm obs}) ] = S^{\rm obs}_N,
\end{equation}
where, we recall, 
$$
\VarR[ K^\star_N(\cdot,\theta_{\rm obs}) ]
=
\frac{\Var[ K^\star_N(\cdot,\theta_{\rm obs}) ]}{\left( \E[ K^\star_N(\cdot,\theta_{\rm obs}) ] \right)^2}.
$$
Given $N$, $K^{\star,{\rm obs}}_N$ and $S^{\rm obs}_N$, our aim is to recover (an approximation of) $\theta_{\rm obs}$. To that aim, we consider the function 
\begin{equation}
\label{eq:def_FNM}
F_{N,M} : \qquad \left\{\begin{aligned}
&  (\xR_+^\star)^2 & \to & \qquad \qquad \qquad \qquad \xR_+   \\
& \theta  & \mapsto   &  \left( \frac{\overline{K}^\star_{N,M}(\theta)}{K^{\star,{\rm obs}}_N}-1 \right)^2 + \left( \frac{S_{N,M}(\theta)}{S^{\rm obs}_N} - 1 \right)^2
\end{aligned}\right.
\end{equation}
which penalizes the sum of the (relative) errors between
\begin{itemize}
\item on the one hand, $\overline{K}^\star_{N,M}(\theta)$ (which is an empirical estimator of $\E \left[ K^\star_N(\cdot,\theta) \right]$ when $M$ is large, see~\eqref{eq:empirical_mean} and~\eqref{eq:def_kstar_N}) and $K^{\star,{\rm obs}}_N$
\item and, on the other hand, $S_{N,M}(\theta)$ (which is an empirical estimator of the relative variance of $K^\star_N(\omega,\theta)$ when $M$ is large, see~\eqref{eq:empirical_var} and~\eqref{eq:def_kstar_N}) and $S^{\rm obs}_N$.
\end{itemize}
Of course, different weights could be assigned to the error on the permeability and the error on its relative variance. We eventually cast our parameter fitting problem in the form of the optimization problem
$$
\inf_{\theta = (\lambda,k) \in (0,\infty) \times {\cal K}} F_{N,M}(\theta),
$$
where ${\cal K} \subset (0,\infty)$ is the admissible set of parameters $k$ such that homogenization holds (even if Assumption~\eqref{eq:ellipticity} is not satisfied for any $k>0$) and the variance of $K^\star_N$ is also well-defined. In the one-dimensional case we focus on in this article, ${\cal K} = (8,\infty)$.

Note that $F_{N,M}(\theta)$ is random, as it depends on the realizations used to evaluate $\overline{K}^\star_{N,M}(\theta)$ and $S_{N,M}(\theta)$ (see~\eqref{eq:empirical_mean} and~\eqref{eq:empirical_var}). For any $\theta$, in the limit when $M \to \infty$, $F_{N,M}(\theta)$ converges almost surely to the deterministic limit
$$
F_N(\theta) = 
\left( \frac{\E \left[ K^\star_N(\cdot,\theta)\right]}{K^{\star,{\rm obs}}_N}-1 \right)^2 
+ 
\left( \frac{\VarR[ K^\star_N(\cdot,\theta) ]}{S^{\rm obs}_N} - 1 \right)^2.
$$
Under Assumption~\eqref{eq:synt}, we have
$$
F_N(\theta) = 
\left( \frac{\E \left[ K^\star_N(\cdot,\theta)\right]}{\E \left[ K^\star_N(\cdot,\theta_{\rm obs})\right]}-1 \right)^2 
+ 
\left( \frac{\VarR[ K^\star_N(\cdot,\theta) ]}{\VarR[ K^\star_N(\cdot,\theta_{\rm obs}) ]} - 1 \right)^2.
$$
When $N \to \infty$, the first term above converges to 
$$
\left( \frac{K^\star(\theta)}{K^\star(\theta_{\rm obs})}-1 \right)^2.
$$ 
For the second term, it is clear that $\VarR[ K^\star_N(\cdot,\theta) ]$ vanishes in the limit $N \to \infty$, since $K^\star_N(\cdot,\theta)$ converges almost surely to a deterministic limit. However, establishing sharp estimates on the rate of convergence is a challenging question. This is why we postpone the discussion on the convergence of the ratio 
$$
\frac{\VarR[ K^\star_N(\cdot,\theta) ]}{\VarR[ K^\star_N(\cdot,\theta_{\rm obs}) ]}
$$ 
to the next section, where we focus on the one-dimensional case, and where more precise results are available. 

\subsection{The one-dimensional case}
\label{sec:inv_1D}

\subsubsection{Theoretical result}
\label{sec:discut}

In the one-dimensional case, we have seen (see~\eqref{eq:astar_1D_weibull}) that
$$
A^\star = \frac{\lambda^4}{\Gamma(1-4/k)}
$$
where $\Gamma$ is the Euler Gamma function. Furthermore, equation~\eqref{eq:approx_hom_1D} implies that
\begin{equation}
\label{eq:rel_var_1D_weibull_pre}
\Var[ A^\star_N ] = \frac{(A^\star)^4}{N} \, \Var \left[ \frac{1}{A(0,\cdot)} \right] + o \left( \frac{1}{N} \right).
\end{equation}
The conductances are distributed according to a Weibull law (see~\eqref{eq:prop}), therefore
$$
\Var[ A^\star_N ] = \frac{\lambda^{16}}{N \, \Gamma(1-4/k)^4} \,  \left( \frac{\Gamma(1-8/k)}{\lambda^8} - \frac{\Gamma(1-4/k)^2}{\lambda^8} \right) + o \left( \frac{1}{N} \right),
$$
hence the relative variance reads
\begin{equation}
\label{eq:rel_var_1D_weibull}
\VarR[ A^\star_N ] = \frac{1}{N} \, \left( \frac{\Gamma(1-8/k)}{\Gamma(1-4/k)^2} - 1 \right) + o \left( \frac{1}{N} \right),
\end{equation}
which implies that
$$
\lim_{N \to \infty}
\frac{\VarR[ A^\star_N(\cdot,\theta) ]}{\VarR[ A^\star_N(\cdot,\theta_{\rm obs}) ]}
=
\frac{
\frac{\Gamma(1-8/k)}{\Gamma(1-4/k)^2} - 1 
}{
\frac{\Gamma(1-8/k_{\rm obs})}{\Gamma(1-4/k_{\rm obs})^2} - 1 
}.
$$
In the one-dimensional case, we are thus able to identify the limit as $N \to \infty$ of $F_N(\theta)$, which reads
\begin{equation}
\label{eq:ideal_functional1D}
F^{\rm 1D}_\infty(\theta) 
:= 
\lim_{N \to \infty} F_N(\theta) 
=
\left( \frac{\lambda^4}{\lambda_{\rm obs}^4}\frac{\Gamma(1 - 4/k_{\rm obs})}{ \Gamma(1-4/k)}-1 \right)^2 
+ 
\left( \frac{
\frac{\Gamma(1-8/k)}{\Gamma(1-4/k)^2} - 1 
}{
\frac{\Gamma(1-8/k_{\rm obs})}{\Gamma(1-4/k_{\rm obs})^2} - 1 
} - 1 \right)^2.
\end{equation}
Obviously, this function is minimal (and vanishes) when $\theta = \theta_{\rm obs}$. It turns out that this minimizer is the unique minimizer, as shown below.

\begin{lmm}
\label{lem:bij}
The function $F_\infty^{\rm 1D}$ defined by~\eqref{eq:ideal_functional1D} has a unique minimizer, which is $\theta_{\rm obs}$.
\end{lmm}

This result is very useful. Homogenization is an averaging process, which filters out many features of the microscopic coefficient $A$. These features cannot be recovered from the knowledge of macroscopic quantities. The above lemma shows (in the one-dimensional case) that, if one assumes a given form for the probability distribution of $A$ (here, a Weibull distribution), then one is able to recover the two parameters of that law on the basis of two macroscopic quantities, the permeability and its relative variance. 

It is also obvious from~\eqref{eq:astar_1D_weibull} that knowing the macroscopic permeability is not enough to uniquely determine the two parameters $\lambda$ and $k$ of the Weibull law. Additional information is needed. Our choice of considering the relative variance of the permeability is motivated by the following observation. This quantity, in the one-dimensional case, only depends (at first order in $N$) on $k$ and does not depend on $\lambda$, as can be seen on~\eqref{eq:rel_var_1D_weibull}. Knowing this quantity is therefore very useful to estimate the parameter $k$. Once $k$ has been identified, knowing the macroscopic permeability yields, using~\eqref{eq:astar_1D_weibull}, an estimation of the parameter $\lambda$.

Of course, it is likely that the knowledge of quantities of interest alternate to the relative variance of the permeability may also prove useful to determine the unknown parameters. Note also that such alternate relevant quantities should be ``different enough'' from the homogenized permeability to indeed bring new information. We do not pursue in that direction. 

\medskip

\begin{figure}[H]
\begin{center}
\includegraphics[scale=0.30]{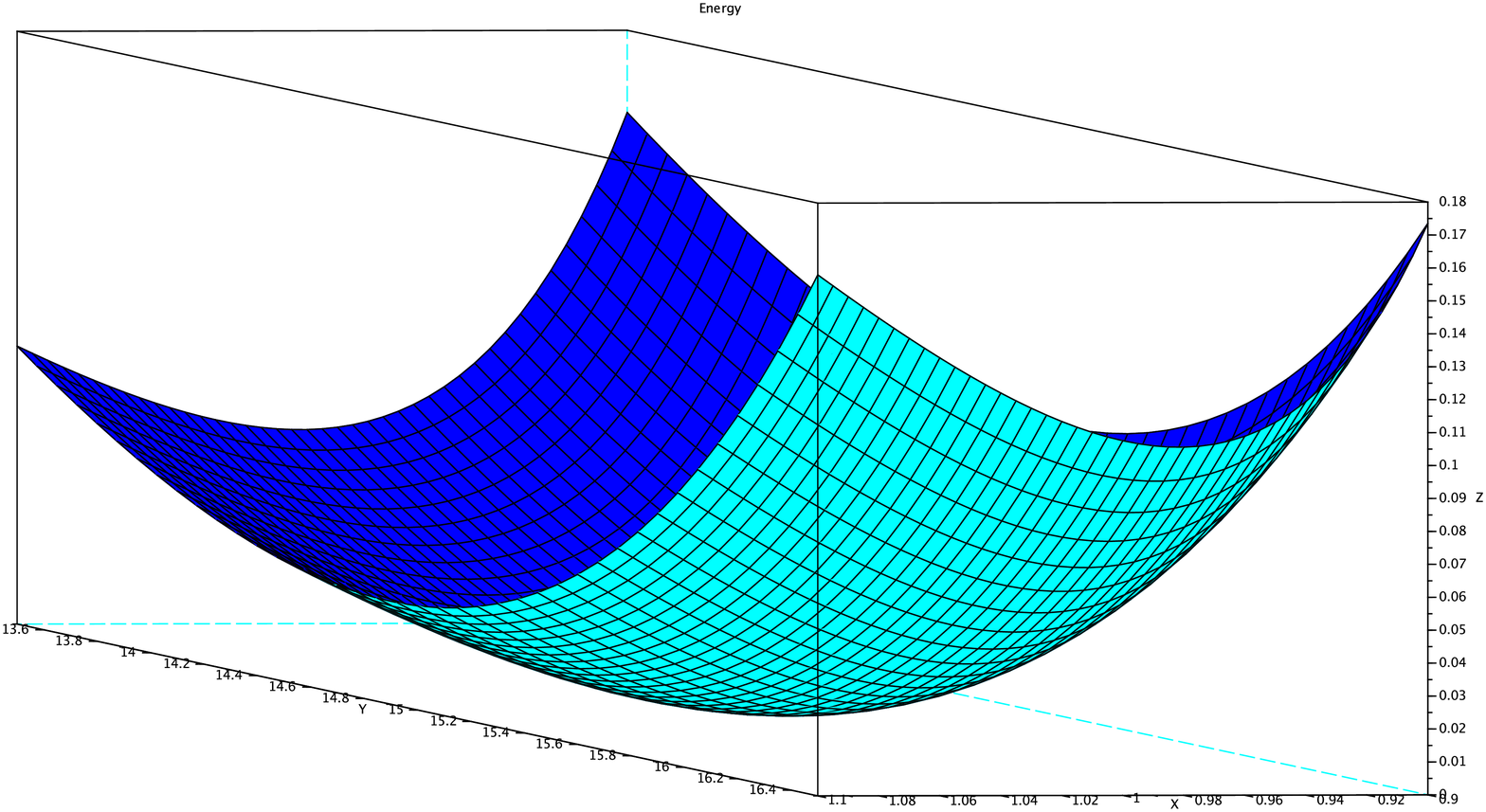}
\end{center}
\caption{Plot of $\theta \mapsto F^{\rm 1D}_\infty(\theta)$ for $\lambda_{\rm obs} = 1$ and $k_{\rm obs} = 15$.}
\label{fig:sampling}
\end{figure}

We plot on Figure~\ref{fig:sampling} the function $\theta \mapsto F^{\rm 1D}_\infty(\theta)$ for $\lambda_{\rm obs} = 1$ and $k_{\rm obs} = 15$. We observe that the function is not degenerated at its minimum, in the sense that its Hessian matrix at $\theta_{\rm obs}$ is positive definite, with eigenvalues equal to 16 and 0.04. We thus expect that a standard algorithm (such as the Newton algorithm) will be able to converge to the minimizer of $F^{\rm 1D}_\infty$. This is indeed the case, as shown in Section~\ref{sec:numerics}.  

\begin{proof}[Proof of Lemma~\ref{lem:bij}.]
The proof consists of three steps: in Step 1, we recall (and prove for the sake of completeness) that $\ln \Gamma$ is a convex function. In Step 2, we prove that the function $$\zeta: k \mapsto \frac{\Gamma(1-8/k)}{\Gamma(1-4/k)^2}$$ is monotone (hence injective). We conclude in Step 3.

\paragraph{\bf Step 1.}
From~\eqref{eq:def_Gamma}, we compute that, for any $z>0$,
\begin{eqnarray*}
\Gamma'(z) &=& \int_0^\infty (\ln t) \, t^{z-1} \, \exp(-t) \, dt,
\\
\Gamma''(z) &=& \int_0^\infty (\ln t)^2 \, t^{z-1} \, \exp(-t) \, dt,  
\end{eqnarray*}
therefore $\Gamma''(z) > 0$ and $\Gamma$ is positive and convex on $(0,\infty)$. In addition, we have
$$
\left( \ln \Gamma \right)''(z) = \frac{\Gamma(z) \Gamma''(z) - (\Gamma'(z))^2}{\Gamma^2(z)} 
$$
which is positive, in view of the Cauchy-Schwartz inequality:
$$
(\Gamma'(z))^2 
= 
\left( \int_0^\infty \left( (\ln t) \, \sqrt{t^{z-1} \, \exp(-t)} \right) \sqrt{t^{z-1} \, \exp(-t)} \, dt \right)^2
<
\Gamma''(z) \Gamma(z).
$$
Therefore, $\ln \Gamma$ is a strictly convex function. 

\medskip

\paragraph{\bf Step 2.}
We define the function
$$
\zeta: k \mapsto \frac{\Gamma(1-8/k)}{\Gamma(1-4/k)^2},
$$
the derivative of which reads
$$
\zeta'(k) = \frac{\Gamma(1-8/k)}{\Gamma(1-4/k)^2} \left(
\frac{8}{k^2}\frac{\Gamma'(1-8/k)}{\Gamma(1-8/k)} 
-
\frac{8}{k^2} \frac{\Gamma'(1-4/k)}{\Gamma(1-4/k)} 
\right).
$$
For any $k>0$, we have that $1 - 8/k < 1 - 4/k$. As a consequence of $\ln \Gamma$ being strictly convex, we have that its derivative is increasing, therefore
$$
\frac{\Gamma'(1 - 8/k)}{\Gamma(1 - 8/k)} < \frac{\Gamma'(1 - 4/k)}{\Gamma(1 - 4/k)}.
$$
We can now conclude that $\zeta'(k) < 0$, hence $\zeta$ is decreasing.

\medskip

\paragraph{\bf Step 3.}
We first observe that
$$
F^{\rm 1D}_\infty(\theta) 
=
\left( \frac{\lambda^4}{\lambda_{\rm obs}^4}\frac{\Gamma(1 - 4/k_{\rm obs})}{ \Gamma(1-4/k)}-1 \right)^2 
+ 
\left( \frac{\zeta(k)-1}{\zeta(k_{\rm obs})-1} - 1 \right)^2.
$$
We obviously have that $\theta_{\rm obs}$ is a minimizer of $F_\infty^{\rm 1D}$, with $F_\infty^{\rm 1D} (\theta_{\rm obs}) = 0$. Conversely, let $\theta$ be a minimizer of $F_\infty^{\rm 1D}$. We thus have $F_\infty^{\rm 1D} (\theta) = 0$, which implies that $\zeta(k) = \zeta(k_{\rm obs})$. The function $\zeta$ being monotone, this implies that $k = k_{\rm obs}$. Since the first term in $F^{\rm 1D}_\infty(\theta)$ also has to vanish, we obtain that $\lambda = \lambda_{\rm obs}$ as well. This concludes the proof.
\end{proof}

\subsubsection{Practical situation}

In the general (i.e. multi-dimensional) case, we have introduced in~\eqref{eq:def_FNM} the function $F_{N,M}$ that we wish to minimize. We have next theoretically identified its limit when $M \to \infty$ and $N \to \infty$. In practice, we cannot take any of these limits, and have thus to work with $F_{N,M}$ defined by
$$
F_{N,M}(\theta) = \left( \frac{\overline{K}^\star_{N,M}(\theta)}{K^{\star, {\rm obs}}_N} -1 \right)^2 + \left( \frac{S_{N,M}(\theta)}{S^{\rm obs}_N} - 1 \right)^2.
$$
In view of~\eqref{eq:rel_var_1D_weibull_pre}, we see that, when $M \to \infty$ and $N \to \infty$, the relative variance $S_{N,M}(\theta)$ is close to 
$$
S_{N,M}(\theta) 
\approx
\frac{(A^\star)^2}{N} \, \Var \left[ \frac{1}{A(0,\cdot)} \right]
=
\frac{\E[W^{-8}]- \left( \E[W^{-4}] \right)^2}{N\,\left( \E[W^{-4}] \right)^2}
=
\frac{\E[W^{-8}]}{N\,\left(\E[W^{-4}]\right)^2} - \frac{1}{N},
$$
where $W$ is a random variable distributed according to the Weibull ${\cal W}(1,k)$. Likewise, 
$$
\overline{K}^\star_{N,M}(\theta) 
\approx
A^\star
=
\frac{\lambda^4}{\E[W^{-4}]}.
$$
Let $\left\{ u_i(\omega) \right\}_{i=1}^N$ be a sequence of i.i.d. random variables uniformly distributed in $[0,1]$. We define 
\begin{equation}
\label{eq:def_wi}
w_i(k,\omega):=(-\ln(1-u_i(\omega)))^{-1/k},
\end{equation}
so that $\left\{ 1/w_i(k,\omega) \right\}_{i=1}^N$ are i.i.d. random variables distributed according to ${\cal W}(1,k)$. In the sequel, we approximate the function to minimize by
\begin{equation}
\label{eq:THE_function}
\widetilde{F}^{\rm 1D}_N(\theta,\omega)
=
\left( \frac{\lambda^4}{K^{\star,{\rm obs}}_N} \left[ \frac{1}{N}\sum_{i=1}^N w_i^4(k,\omega)\right]^{-1}-1 \right)^2 
+ 
\left(\frac{1}{S^{\rm obs}_N} \left[\frac{\sum_{i=1}^N w_i^8(k,\omega)}{\left(\sum_{i=1}^N w_i^4(k,\omega)\right)^2}-\frac{1}{N} \right]-1\right)^2. 
\end{equation}
This function is consistent in the sense that it almost surely converges, when $N \to \infty$, to the exact function~\eqref{eq:ideal_functional1D}. On the other hand, $\widetilde{F}^{\rm 1D}_N(\theta,\omega)$ is random, and thus somewhat mimics the difficulties that one would encounter in the multi-dimensional case when working with $F_{N,M}(\theta)$.

\section{Numerical results} 
\label{sec:numerics}

We briefly explain in Section~\ref{sec:optimization} how in practice we minimize the function~\eqref{eq:THE_function}, before turning in Section~\ref{sec:num} to our numerical results. As pointed out in the introduction, we only consider here the one-dimensional case, and postpone the study of two-dimensional examples to the future work~\cite{LMOSxx}.

\subsection{Optimization algorithm}
\label{sec:optimization}

We show in Appendix \ref{app:calcul} how to compute the first and second derivatives of the function $\widetilde{F}^{\rm 1D}_N(\theta,\omega)$ defined by~\eqref{eq:THE_function} with respect to $\theta = (\lambda,k)$. We are thus in position to use the Newton algorithm, and compute a sequence $\theta_j$ according to
\begin{equation}
\label{eq:newton}
\theta_{j+1} = \theta_j - \mu_j \left[ {\cal H} \left( \widetilde{F}^{\rm 1D}_N \right)(\theta_j) \right]^{-1} \nabla \widetilde{F}^{\rm 1D}_N(\theta_j), 
\end{equation}
where ${\cal H} \left( \widetilde{F}^{\rm 1D}_N \right) \in \xR^{2 \times 2}$ is the Hessian matrix of $\widetilde{F}^{\rm 1D}_N$ and $\nabla \widetilde{F}^{\rm 1D}_N \in \xR^2$ is the gradient of $\widetilde{F}^{\rm 1D}_N$ (for the sake of simplicity, we keep implicit the dependence with respect to $\omega$). In turn, $\mu_j > 0$ is the step-size by which we move. To choose $\mu_j$, we have used a line-search algorithm (along the descent direction prescribed by the Newton algorithm) using Goldstein (respectively Armijo) rule to increase (respectively decrease) the step-size.

\medskip

We note that the function $\theta \mapsto F^{\rm 1D}_\infty(\theta)$ is not convex. It is possible to find some $\theta$ such that the Hessian matrix ${\cal H} \left( F^{\rm 1D}_\infty \right)(\theta)$ is not positive definite, but rather has (at least) one negative eigenvalue. We thus cannot expect the function $\theta \mapsto \widetilde{F}^{\rm 1D}_N(\theta)$ to be convex (even for large values of $N$), and the Newton algorithm to be globally convergent. We are therefore careful to start the Newton iterations from an initial guess $\theta_0$ (given by physical experiments) that we hope to be close enough to the minimizer of $\widetilde{F}^{\rm 1D}_N$.

\subsection{Numerical results}
\label{sec:num}

In all what follows, we set $N=10^5$.

\subsubsection{Robustness of the algorithm with respect to the initial guess}

Our first numerical test is a simple one, to check whether the Newton algorithm~\eqref{eq:newton} is indeed able to minimize the function $\theta \mapsto \widetilde{F}^{\rm 1D}_N(\theta,\omega)$. We pick once for all one realization of the i.i.d. random variables $\left\{ u_i(\omega) \right\}_{1 \leq i \leq N}$ (which, we recall, are uniformly distributed in $[0,1]$). We then build $\left\{ w_i(k,\omega) \right\}_{1 \leq i \leq N}$ according to~\eqref{eq:def_wi} and consider the function $\theta \mapsto \widetilde{F}^{\rm 1D}_N(\theta,\omega)$ defined by~\eqref{eq:THE_function}, where the observed quantities are defined by
$$
K^{\star,{\rm obs}}_N
=
\lambda^4_{\rm obs} \left[ \frac{1}{N}\sum_{i=1}^N w_i^4(k_{\rm obs},\omega)\right]^{-1},
\qquad
S^{\rm obs}_N
=
\frac{\sum_{i=1}^N w_i^8(k_{\rm obs},\omega)}{\left(\sum_{i=1}^N w_i^4(k_{\rm obs},\omega)\right)^2}-\frac{1}{N},
$$
with $\lambda_{\rm obs} = 1$ and $k_{\rm obs} = 15$. The function $\theta \mapsto \widetilde{F}^{\rm 1D}_N(\theta,\omega)$ obviously vanishes at $\theta_{\rm obs} = (\lambda_{\rm obs},k_{\rm obs})$. 

\medskip

We run the Newton algorithm~\eqref{eq:newton} starting from several initial guesses $\theta_0$, and check that it indeed always converges to $\theta_{\rm obs}$ in a limited number of iterations. We also observe that, for some initial guesses, using an adaptive step-size $\mu_j$ as in~\eqref{eq:newton} is critical: if, in contrast, one uses the step-size $\mu_j=1$, then the algorithm may not converge, or converges within a much larger number of iterations.

\subsubsection{Robustness with respect to statistical noise}

For our second test, we proceed as follows. We first set $\theta_{\rm ref} = (\lambda_{\rm ref},k_{\rm ref}) = (1,15)$ and pick one realization of the i.i.d. random variables $\left\{ u_i(\overline{\omega}) \right\}_{1 \leq i \leq N}$ (which, we recall, are uniformly distributed in $[0,1]$). We then build $\left\{ w_i(k_{\rm ref},\overline{\omega}) \right\}_{1 \leq i \leq N}$ according to~\eqref{eq:def_wi} and define once for all the macroscopic observed quantities as
\begin{equation}
\label{eq:synt2}
K^{\star,{\rm obs}}_N
=
\lambda^4_{\rm ref} \left[ \frac{1}{N}\sum_{i=1}^N w_i^4(k_{\rm ref},\overline{\omega})\right]^{-1},
\qquad
S^{\rm obs}_N
=
\frac{\sum_{i=1}^N w_i^8(k_{\rm ref},\overline{\omega})}{\left(\sum_{i=1}^N w_i^4(k_{\rm ref},\overline{\omega})\right)^2}-\frac{1}{N}.
\end{equation}
We now fix the initial guess $\theta_0=(1.1,16.5)$ (10\% off the reference value $\theta_{\rm ref}$) and set $M=500$. For any $1 \leq m \leq M$, we perform the following procedure:
\begin{itemize}
\item we draw a realization of $N$ i.i.d. random variables $\left\{ u_i(\omega_m) \right\}_{1 \leq i \leq N}$ which is independent of the realization $\left\{ u_i(\omega_{m'}) \right\}_{1 \leq i \leq N}$ for any $m' \neq m$, and independent of the realization $\left\{ u_i(\overline{\omega}) \right\}_{1 \leq i \leq N}$ used to compute $K^{\star,{\rm obs}}_N$ and $S^{\rm obs}_N$ in~\eqref{eq:synt2};
\item using $\left\{ u_i(\omega_m) \right\}_{1 \leq i \leq N}$, we build $w_i(k,\omega_m)$ according to~\eqref{eq:def_wi} and we consider the function $\theta \mapsto \widetilde{F}^{\rm 1D}_N(\theta,\omega_m)$ defined by~\eqref{eq:THE_function}, i.e. 
$$
\widetilde{F}^{\rm 1D}_N(\theta,\omega_m)
=
\left( \frac{\lambda^4}{K^{\star,{\rm obs}}_N} \left[ \frac{1}{N}\sum_{i=1}^N w_i^4(k,\omega_m)\right]^{-1}-1 \right)^2 
+ 
\left(\frac{1}{S^{\rm obs}_N} \left[\frac{\sum_{i=1}^N w_i^8(k,\omega_m)}{\left(\sum_{i=1}^N w_i^4(k,\omega_m)\right)^2}-\frac{1}{N} \right]-1\right)^2. 
$$
Recall that the macroscopic observed quantities are independent of $\omega_m$.
\item we run the Newton algorithm~\eqref{eq:newton} to minimize the function $\theta \mapsto \widetilde{F}^{\rm 1D}_N(\theta,\omega_m)$. The optimal parameter found by the algorithm depends on $\omega_m$ and is denoted $\theta_{\rm opt}(\omega_m)$. Since the realization $\omega_m$ is different from the reference realization $\overline{\omega}$, we have in general $\theta_{\rm opt}(\omega_m) \neq \theta_{\rm ref}$. 
\end{itemize}

\begin{figure}[H]
\includegraphics[scale=0.4]{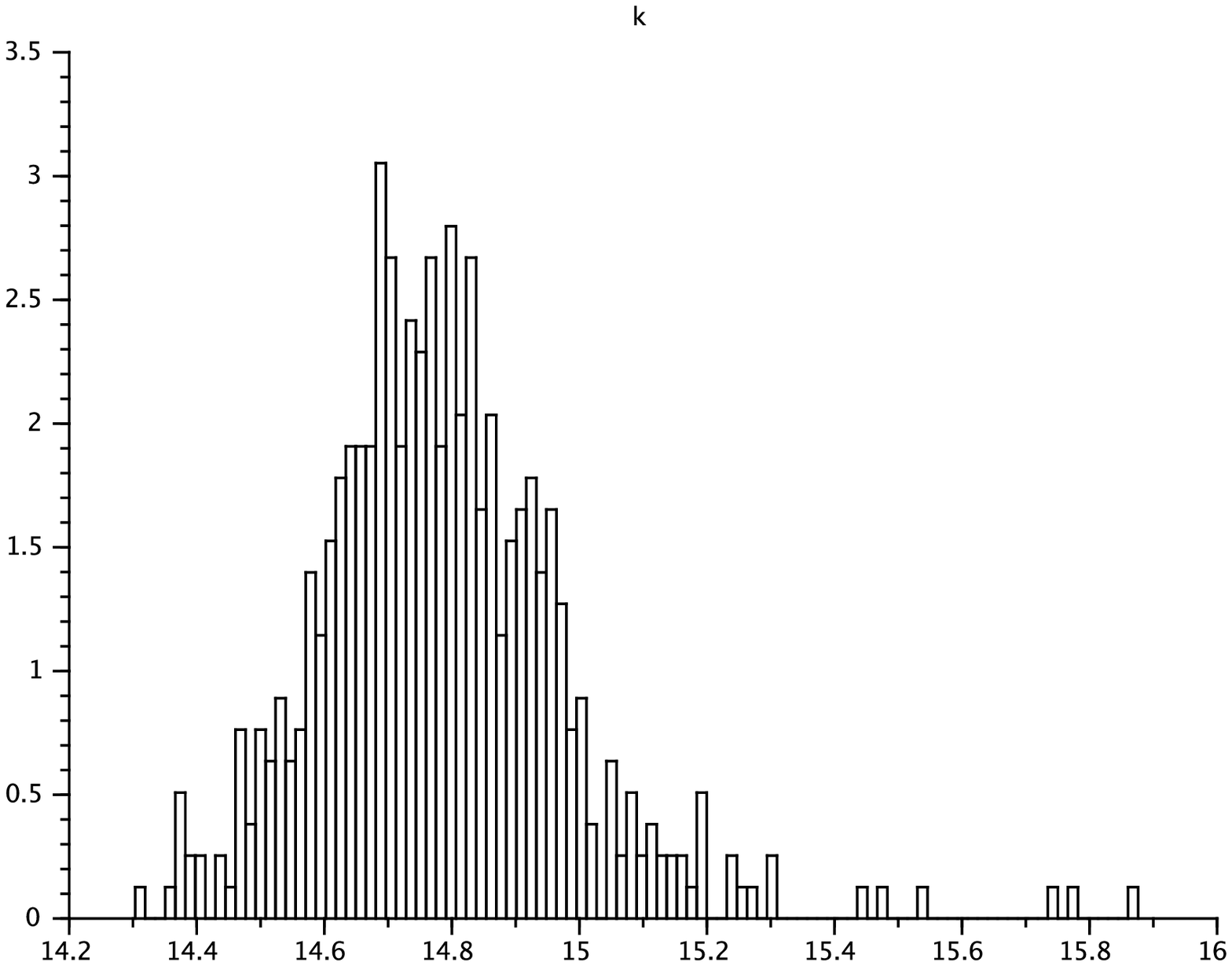}
\includegraphics[scale=0.4]{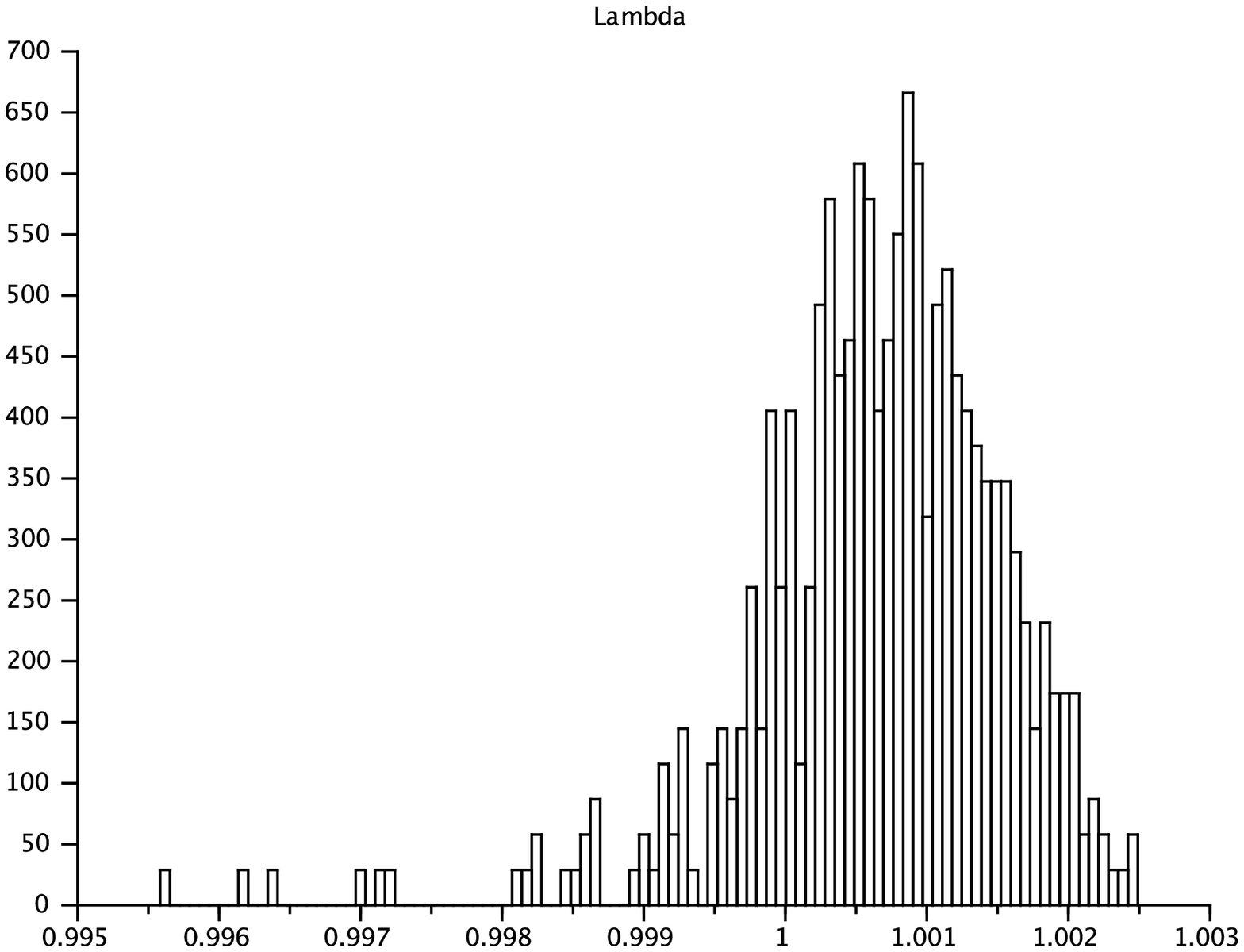}
\caption{Top: distribution of $k_{\rm opt}(\omega)$. Bottom: distribution of $\lambda_{\rm opt}(\omega)$.}
\label{fig:distribution_theta}
\end{figure}

We show on Figure~\ref{fig:distribution_theta} the histogram of the optimal parameters $\theta_{\rm opt}(\omega_m)$ for $1 \leq m \leq M$. We see that these histograms are centered close to the reference value ($k_{\rm ref}$, resp. $\lambda_{\rm ref}$). There is however a small bias, i.e. $\E\left( \theta_{\rm opt} \right) \neq \theta_{\rm ref}$. We also observe that the width of these histograms (related to the variance of $k_{\rm opt}$ and $\lambda_{\rm opt}$) is quite small. 

\begin{rmrk}
Of course, the variance of $k_{\rm opt}$ and $\lambda_{\rm opt}$ is related to $N$. In the limit $N \to \infty$, the function $\widetilde{F}^{\rm 1D}_N(\theta,\omega)$ almost surely converges to a deterministic limit, and we thus expect $k_{\rm opt}$ and $\lambda_{\rm opt}$ to almost surely converge to a deterministic limit. But this is not the regime we are interested in, since in practice (in the two-dimensional case), we have to work with the {\em random} function $F_{N,M}$.
\end{rmrk}

\begin{figure}[H]
\includegraphics[scale=0.4]{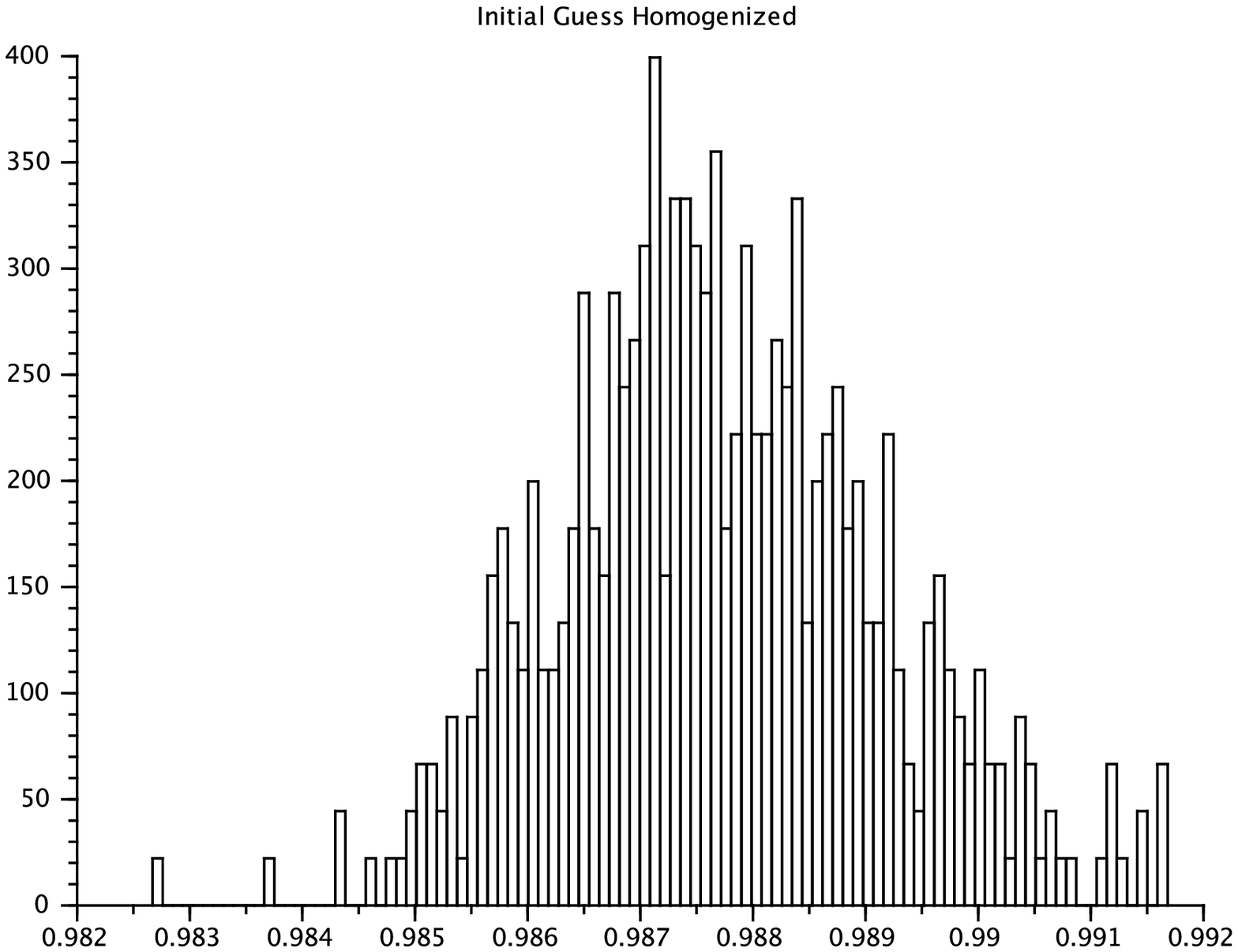}
\includegraphics[scale=0.4]{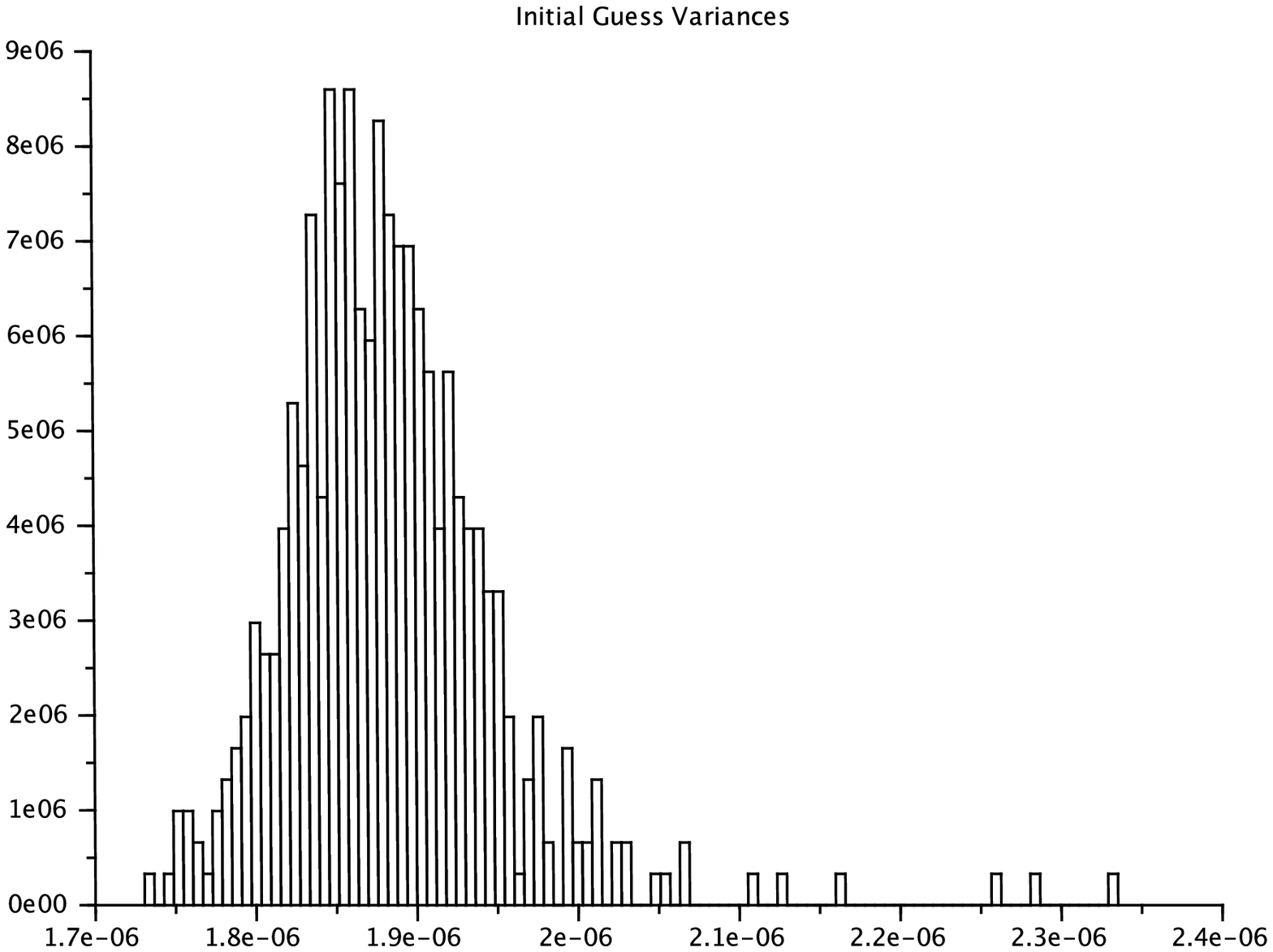}
\caption{Top: distribution of $K^\star_N(\theta_0,\omega)$. Bottom: distribution of $S_N(k_0,\omega)$.}
\label{fig:distribution_IG}
\end{figure}

We next compare the variance of $\theta_{\rm opt}$ with the amount of randomness introduced in the function $\widetilde{F}^{\rm 1D}_N(\cdot,\omega)$ defined by~\eqref{eq:THE_function}. By construction,
$$
\widetilde{F}^{\rm 1D}_N(\theta,\omega)
=
\left( \frac{K^\star_N(\theta,\omega)}{K^{\star,{\rm obs}}_N}-1 \right)^2 
+ 
\left(\frac{S_N(k,\omega)}{S^{\rm obs}_N} -1\right)^2
$$
with
$$
S_N(k,\omega)
=
\left[\frac{\sum_{i=1}^N w_i^8(k,\omega)}{\left(\sum_{i=1}^N w_i^4(k,\omega)\right)^2}-\frac{1}{N} \right],
$$
which is an approximation of the relative variance of $K^\star_N(\theta,\omega)$. We show on Figure~\ref{fig:distribution_IG} the histograms, for $1 \leq m \leq M$, of $K^\star_N(\theta_0,\omega_m)$ and of $S_N(k_0,\omega_m)$, for the initial guess parameter $\theta_0=(1.1,16.5)$. 

On this test-case, we compute that $\Var[\lambda_{\rm opt}] \approx 7.9 \ 10^{-7}$ and $\Var[k_{\rm opt}] \approx 3.8 \ 10^{-2}$, thus
$$
\VarR[\lambda_{\rm opt}] \approx 7.9 \ 10^{-7} \quad \text{ and } \quad \VarR[k_{\rm opt}] \approx 1.7 \ 10^{-4}.
$$
On the other hand, $A^\star(\theta_0) \approx 1.2$, $\Var \left[ A^\star_N(\theta_0) \right] \approx 2.0 \ 10^{-6}$ and $\Var \left[ S_N(k_0) \right] \approx 4.5 \ 10^{-15}$, thus
$$
\VarR \left[ K^\star_N(\theta_0) \right] \approx 1.4 \ 10^{-6} \quad \text{ and } \quad \VarR \left[ S_N(k_0) \right] \approx 10^{-3}.
$$
We thus observe that the relative variance of the optimal parameters is roughly of the same order of magnitude as the relative variance introduced in the function to minimize. Given the amount of noise present in the system, our procedure robustly identifies the optimal parameters of the microscopic distribution.

\appendix

\section{Computation of the derivatives of~\eqref{eq:THE_function}}
\label{app:calcul}

We introduce 
$$
f(\lambda,k) := \lambda^4\left(\sum_{i=1}^N w_i^4(k)\right)^{-1}
$$
and 
$$
g(k) := \left(\sum_{i=1}^N w_i^8(k)\right)\left(\sum_{i=1}^Nw_i^4(k)\right)^{-2}
$$
where $w_i(k)$ is defined by~\eqref{eq:def_wi}, and recast the function~\eqref{eq:THE_function} as
$$
\widetilde{F}^{\rm 1D}_N(\theta)
=
\left(\frac{N}{K^{\star,{\rm obs}}_N} \ f(\lambda,k)-1\right)^2
+
\left(\frac{1}{S^{\rm obs}_N}\left(g(k)-\frac{1}{N}\right)-1\right)^2,
$$
where we have kept implicit the dependence with respect to $\omega$. Computing the derivatives of $\widetilde{F}^{\rm 1D}_N$ therefore amounts to computing those of $f$ and $g$. 

\medskip

A tedious but straightforward computation leads to the following expressions:
\begin{align*}
g'(k)
&=
\frac{8}{k}\left(\sum_i w_i^4\right)^{-2}\left[\frac{\sum_i w_i^8}{\sum_i w_i^4} \sum_i \ln(w_i)w_i^4-\sum_i\ln(w_i)w_i^8\right], 
\\
g''(k)
&=
\frac{16}{k^2}\left(\sum_iw_i^4\right)^{-2}\left[\sum_i \ln(w_i)w_i^8(1+4\ln w_i)- \frac{\sum_i w_i^8}{\sum_i w_i^4} \sum_i\ln(w_i)w_i^4(1+2\ln w_i)\right] 
\\
& \qquad 
-\frac{32}{k^2}\left(\sum_i w_i^4\right)^{-3}\left[ 4\left(\sum_i w_i^8\ln(w_i)\right)\left(\sum_i w_i^4\ln(w_i)\right)-3\left(\sum_i w_i^4\ln(w_i)\right)^2\frac{\sum_i w_i^8}{\sum_i w_i^4} \right],
\end{align*}
whereas
\begin{align*}
\partial_\lambda f
& = 
\frac{4}{\lambda}f(\lambda,k), 
\\
\partial^2_{\lambda\lambda}f 
& =
\frac{12}{\lambda^2}f(\lambda,k),
\\
\partial_k f
& = 
\frac{4\lambda^4}{k}\left(\sum_i \ln(w_i)w_i^4\right) \left(\sum_i w_i^4\right)^{-2},
\\
\partial^2_{k k}f 
& = 
-\frac{8\lambda^4}{k^2}\left(\sum_i w_i^4\right)^{-2}\left[\sum_i\ln(w_i)w_i^4(1+2\ln w_i)\right]+\frac{32\lambda^4}{k^2}\left(\sum_iw_i^4\right)^{-3}\left(\sum_i\ln(w_i)w_i^4\right)^2, 
\\
\partial^2_{\lambda k} f
& =  
\frac{16\lambda^3}{k}\left(\sum_i \ln(w_i)w_i^4\right) \left(\sum_i w_i^4\right)^{-2}=\frac{4}{\lambda}\partial_k f.
\end{align*}

\medskip

\thanks{\textit{Acknowledgements.} We thank Tony Leli\`evre for giving us the opportunity to work on this problem during the CEMRACS 2013 ({\tt http://smai.emath.fr/cemracs/cemracs13/}), and Nicolas Champagnat, Tony Leli\`evre and Anthony Nouy for the organization of this event. We are grateful to Samir B\'ekri, Daniel Coelho, Claude Le Bris, Tony Leli\`evre and Benjamin Rotenberg for enlightening discussions. 
The work of FL and WM is partially supported by ONR under Grant 
N00014-12-1-0383. WM gratefully acknowledges the support from Labex MMCD
(Multi-Scale Modelling \& Experimentation of Materials for Sustainable
Construction) under contract ANR-11-LABX-0022. WM, AO and MS acknowledge financial support from NEEDS {\em Milieux Poreux}.}

\end{document}